\documentclass{article}

                                                                                                                                                                                                                                                                                             \usepackage{verbatim} %Allows use of \begin{comment} ... \end{comment}
\usepackage{shapepar,bm}
\usepackage{setspace}
\usepackage[dvips,arrow,matrix,ps,color,line]{xy}
\usepackage{theorem,amsfonts,amssymb,amscd, amsmath}
\usepackage{geometry}
\usepackage{makeidx}

\usepackage{courier}
\usepackage{mathpazo}			% tt
\usepackage{eulervm}			% a better implementation of the euler package (not in gwTeX)
\usepackage[bb=boondox,bbscaled=1.05,scr=dutchcal]{mathalfa}	%
\usepackage{dsfont}
\normalfont
\usepackage{amsmath}
\usepackage{stmaryrd}
\newcounter{paranum}[section]{}

\usepackage{stmaryrd}  % for  \llbracket   \rrbracket \llparenthesis \rrparenthesis   \llceil \rrceil         \llfloor \rrfloor  

\usepackage{graphicx}
\usepackage{graphics}
\usepackage{tikz-cd}
\usepackage{epstopdf}

\usepackage{amssymb,graphics,hyperref}

\hypersetup{
    colorlinks = true,
    linkcolor = blue,
    anchorcolor = red,
   citecolor = blue,
    filecolor = red,
    pagecolor = red,
    urlcolor = blue}

\newcommand{\mathsym}[1]{{}}

\makeindex

\numberwithin{equation}{section}

\def\beqns{\begin{eqnarray*}}
\def\eqns{\end{eqnarray*}}
\def\mScal{\mathbb S}
\newcommand{\llb}{\llbracket}
\newcommand{\rrb}{\rrbracket}

\def\bmtx{\begin{matrix}}
\def\emtx{\end{matrix}}

\def\bd{\bm \d}

\def\DD{{\mathbb D}}
\def\H24bft2{H_{2,4}(\bft_2)}

\def\NN{\mathbb N}
\def\bqna{\begin{eqnarray*}}
\def\eqna{\end{eqnarray*}}
\def\sI{t_1+t_2}
\def\sD{t_1^2+t_1t_2+t_2^2}
\def\sII{t_1t_2}
\def\sDI{t_1^2t_2+t_1t_2^2}
\def\sDD{t_1^2t_2^2}

\def\DJKM{Date, Jimbo, Kashiwara and Miwa}

\def\bfz{{\mathbf z}}

\def\ovsig{\overline{\sigma}}

\def\bfX{{\mathbf X}}

\def\bfc{{\mathbf c}}
\def\mc{c}
\def\bfmc{{\mathbf c}}

\def\MR{MR}
\def\mS{{ S}}
\def\bfmS{{\mathbf S}}

\def\tr{\mathrm{tr}}

\def\bfv{{\mathbf v}}
\def\bfw{{\mathbf w}}

\def\Ccal{{\mathcal C}}

\def\d{\partial}

\def\bfx{{\mathbf x}}

\def\sgn{\mathrm{sgn}}

\def\ZZ{\mathbb Z}

\def\CC{\mathbb C}

\def\Ecal{{\mathcal E}}

\def\QQ{\mathbb Q}
\def\PP{\mathbb P}

\def\cocoa{{\hbox{\rm C\kern-.13em o\kern-.07em C\kern-.13em o\kern-.15em A}}}

\def\GG{\mathbb G}

\def\GG{{\bf G}}

\def\bft{{\bf t}}

\def\End{\mathrm{End}}

\def\blamb{{\bm \lambda}}
\def\bmu{{\bm \mu}}

\def\Pcal{{\mathcal P}}

\def\w2M{\bigwedge^2M}

\def\w{\wedge }
\def\bw{\bigwedge }

\protect
\protect
\protect
\protect

\protect
\protect

\def\sra{\rightarrow}
\def\lra{\longrightarrow}

\def\proof{\noindent{\bf Proof.}\,\,}
\def\qed{{\hfill\vrule height4pt width4pt depth0pt}\medskip}
\def\be{\begin{equation}}
\def\ee{\end{equation}}
\def\bclm{\begin{claim}}
\def\eclm{\end{claim}}
\def\beqn{\begin{eqnarray}}
\def\eeqn{\end{eqnarray}}
\def\beqn*{\begin{eqnarray*}}
\def\eeqn*{\end{eqnarray*}}

\def\red{\textcolor[rgb]{1.00,0.00,0.00}}
\theoremstyle{change}
\theorembodyfont{\rmfamily}

\newtheorem{claim}{}[section]
\global
\global

\def\no@breaks#1{{\def\\{ \ignorespaces}#1}}    % disallow explicit line breaks

  \makeatletter
\def\cleardoublepage{\clearpage\if@twoside \ifodd\c@page\else
\hbox{} \thispagestyle{empty}
\newpage
\if@twocolumn\hbox{}\newpage\fi\fi\fi} \makeatother

\usepackage{eso-pic,graphicx}
\makeatletter
\newcommand\BackgroundPicture[2]{%
  \setlength{\unitlength}{1pt}%
  default \put(0,\strip@pt\paperheight){%
  \parbox[t][\paperheight]{\paperwidth}{%
    \vfill
     \centering \includegraphics[angle=#2, width=15cm, height=15cm,  bb=0 0 150 150]{#1}
    \vfill
}}} %
\makeatother
\usepackage{palatino}

\title{Action of Free Fermions on Symmetric Functions }

\author{\sc{Letterio Gatto} \footnote{Dipartimento di Scienze Matematiche, Politecnico di Torino, Corso Duca degli Abruzzi 24, 10129, Torino (TO)} \,\,\& \sc{Malihe Yousofzadeh} \footnote{
{Department of Pure Mathematics, Faculty of Mathematics and Statistics, University of Isfahan, Isfahan, P.O. Box 81746-73441, Iran and
School of Mathematics, Institute for Research in Fundamental Sciences (IPM), P.O. Box: 19395-5746, Tehran, Iran}}}

\date{}

\begin{document}

\maketitle
\noindent
{\bf Abstract.} The Clifford algebra of the endomorphisms of the exterior algebra of a countably dimensional vector space induces natural {\em bosoni}c shadows, i.e. families of  linear maps between the cohomologies of complex Grassmannians. The main result of this paper is to provide a determinantal formula expressing generating functions of such endomorphisms unifying  several classical special cases. For example the action over a point recovers the Jacobi-Trudy formula in the theory of symmetric functions or the Giambelli's one in classical Schubert calculus, whereas the action of degree preserving endomorphisms take into account a finite type version of the Date-Jimbo-Kashiwara-Miwa bosonic vertex operator representation of the Lie algebra $gl(\infty)$. The fermionic actions on (finite type) bosonic spaces is described in terms of the classical theory of symmetric functions. The main guiding principle is the fact that the exterior algebra is a (non irreducible) representation of the ring of symmetric functions, which is the  way we use  to spell the ``finite type'' Boson-Fermion correspondence.

\tableofcontents
%vectors.  

\newpage
%\doublespacing
\setstretch{1.80}
\section{Introduction} 
\claim{} Let $r,n\in (\NN\cup\{\infty\})^2$ such that $0\leq r\leq n$. The main purpose of this paper is to describe linear maps between singular cohomologies 
$B_{r,n}$  of complex Grassmannians 
$G(r,n)$ induced by homogeneous endomorphisms (basically the {\em free fermions} that the title alludes to) of the exterior algebra $\bw\QQ^n$. By {\em free fermions} one usually means the elements of the Clifford algebra $
\Ccal_\infty$ 
spanned by $1$ and $(\psi_i, \psi_i^*)_{i\in\ZZ}$, subjects to the relations $\psi_i\psi_j+
\psi_i\psi_j=0$, $\psi^*_i\psi^*_j+\psi^*_i\psi^*_j=0$ and $\psi_i\psi_j^*+
\psi_j^*\psi_i=\delta_{ij}$, see e.g. \cite{jimbomiwa} or \cite[Section 1]{SavageBF}. The {\em elementary endomorphisms} of the 
exterior algebra $\bw \QQ^n$ ($n\in\NN\cup\{\infty\}$), when multiplied, are subject to the same commutation relations of the Clifford algebra $
\Ccal_\infty$ alluded to above. By wedging and contracting, they define linear maps $
\bw^r\QQ^n\sra \bw^{r+s}\QQ^n$  and, because of a certain {\em finite type boson-fermion correspondence}, linear 
maps $B_{r}\sra B_{r+s}$,  whence, up to a slight  abuse 
of terminology, the title of the paper.  The main motivation to pursue this  goal is that it place itself at the crossroad of a number of subjects, such as the theory of symmetric functions, Schubert Calculus for Grassmannians, as well as the representation theory of Lie algebras of endomorphisms of infinite dimensional vector space à la \DJKM (DJKM) \cite{DJKM01}. 

The simplest cases are those corresponding to the pair $(1,n)$, for any $n$,  and $(r,\infty)$, for any $r$, 
which will be partly discussed separately via ad hoc arguments in section 
\ref{sec:specialc}. The case $(\infty,\infty)$ is related to the DJKM bosonic vertex operator representation of the Lie algebra  $gl_\infty(\QQ)$ of  matrices $A:
\ZZ\times \ZZ\sra \QQ$, with finitely many non-zero diagonals.

The general instance,  especially the trickiest one for finite $r$ and $n$, is taken into account by our main 
result, Theorem \ref{thm4:mainthm}, built out of  a new determinantal expression which on one hand  can be seen as a 
generalization of the classical Giambelli's formula of Schubert calculus and, on the other,  keeps 
into account, as a particular case, the structure of the $\QQ$--polynomial ring  in $r$ 
indeterminates as a representation of the Lie algebra  $gl_\infty(\QQ)$. To  achieve the latter one 
may exploit suitable operators, carrying a classical flavor,
whose limit for $r\sra \infty$ are precisely those used  by \DJKM\, (DJKM), in 
\cite{DJKM01}, to 
obtain the bosonic vertex operator representation of the Lie algebra $gl_\infty(\CC)$ 
(see also \cite{jimbomiwa}).  

To this respect, it should be mentioned that Laksov and Thorup, in their reworking of Schubert calculus in terms of action on an exterior power,  in the paper \cite{LakTh01} affirm that ``in the work of E. Date, M.~Jimbo, M.~Kashiwara, and T.~Miwa [4] Schur fucnctions appear in connection with exterior power in another context''. The present paper, on the other hand, is a further indication, besides e.g. \cite{BeNa, BeGa} that the context is indeed  the same and that the celebrated {\em bosonic vertex operator representation}of \DJKM\,\, could  be seen as an instance of the general Schubert Calculus studied in \cite{LakTh01} living on the Universal Grassmann Manifold (inductive limit of finite dimensional grassmannians) by Sato \cite{sato1}

The underlying conceptual framework  is easy. If  $r\leq 
n <\infty$, the $gl_n(\QQ)$-module structure of $B_{r,n}$ is induced by the standard 
representation of the $n\times n$ square matrices on the exterior power $\bw^r\QQ^n$, due 
to the linear isomorphism $B_{r,n}\sra \bw^r\QQ^n$ . The latter can be understood, as we 
will, as a sort of  ``finite 
dimensional'' version of the {\em boson-fermion} correspondence, for which there are now 
many references in the literature { (see e.g. the relatively recent papers \cite{Jing,Rosas})}, and which reflects the fact that  
both spaces possess bases of the same cardinality. 
%This carries  a very classical flavor, as  it may be 
%understood in terms of the well known isomorphism mapping all the symmetric polynomials 
%to skew-symmetric ones via multiplication with the $r\times r$ Vandermonde determinant, 
%playing the role of the {\em vacum 
%%vector}. 
For finite $r$ and $n$ (in contrast to 
$r=n=\infty$ as in the DJKM picture), the determination of an 
explicit generating function of elementary matrices is tricky. One way to cope with the 
issues is proposed in \cite{gln}. On the other hand one can  notice  that 
$B_{r,n}$ is more generally a  module over the algebra of the degree $0$ homogeneous  
endomorphisms of $\bw\QQ^n$, a Lie subalgebra of $\End_\QQ(\bw\QQ^n)$. This was studied 
in \cite{BeCoGaVi} for $n=\infty$, generalizing formulas deduced in \cite{gln}. The structure of 
$B=B_{\infty,\infty}$ as a representation of the endomorphisms of the exterior algebra $
\bw\QQ^\infty$ is considered in \cite{BeGa}: the output are formulas involving  a vertex 
operators description of the generating function which reduces to the DJKM ones  when the exterior algebra endomorphisms are induced by
$gl_\infty(\QQ)$ via the {\em trace representation} (see Sec. \ref{sec3:secsig}). Exploiting the 
remarkable fact that the exterior algebra of a countably infinite vector space is in turn
a representation of the ring of symmetric functions, the aim of this   paper is  
to provide a unified framework to all the previously mentioned situations. 
The finite-dimensional boson-fermion correspondence yields  the general determinantal formula which Theorem \ref{thm4:mainthm} accounts for.

\claim{}  From a representation-theoretic perspective, the boson-fermion correspondence -- first introduced in quantum field theory (see, e.g., \cite{FrenkQFT})-- is merely an isomorphism between two modules over the infinite-dimensional Heisenberg algebra (see e.g. \cite[Section~1]{SavageBF} for a concise enlightening overview). The former, the bosonic 
Fock space, is  $B[q, q^{-1}]$, where  $B$ is the free polynomial algebra generated by 
infinitely many indeterminates $(x_1,x_2,\ldots)$ (the bosonic Fock space in charge $0$), whereas the latter, the {\em fermionic}  one,  is usually  understood as the linear span of the so-called charged  partitions, which in turns can be thought of as 
monomials of a semi-infinite exterior 
power of an infinite dimensional vector space  $V$ (see \cite{Tingley, KacPet, LicAsA} and also \cite[Section~4.1]{SDIWP}). In charge zero, the 
boson-fermion correspondence maps a basis element of $B$,  parametrized by a partition, to an 
infinite exterior monomial attached to the same partition (see e.g. \cite[Theorem 6.1]
{KacRaRoz}). The irreducible $B$-module structure of the fermionic space allows to deal with a wealth of 
related 
situations in mathematical physics and in (algebraic) geometry. Geometric realization of the fermionic Fock spaces and Clifford and Heisenberg operators, after Nakajima and Grojnowskij, are faced e.g. in \cite{SavageBF}, also in connection with the geometry of moduli spaces of framed torsion-free sheaves as in \cite{LicAsA}. From the mathematical physics point of view, one must mention soliton equations,  the celebrated Kadomtsev-Petshiasvili 
(KP) hierarchy, the bosonic vertex operators representation of infinite dimensional Lie 
algebras \cite{DJKM01,jimbomiwa} {  and Jing's article \cite{Jing}}, not to speak about the representation theory of Kac--Moody algebras.  In particular,
the KP-hierarchy rules 
equations of the group orbit of the highest weight vector of $B$ into the projective space 
of the corresponding highest weight module;  in addition, the fact that the polynomial 
solutions of the 
celebrated    KP-hierarchy are parametrized by the points of the Sato's Universal Grassmann 
Manifold (UGM), as constructed e.g. in \cite{sato1}, can be recovered via the bosonic 
representation  of 
the Clifford algebra $\Ccal_\infty$ of free fermions; see  
\cite{dimitrovGran,EswaFut01,Eswa01,Fut01,KaWaki01} and \cite{Yousof08,Yousof09,Yousof10}  
for further information on highest weight modules and 
other kind of representations of affine  Kac-Moody Lie (super)algebras.
Most 
of the above situations have been considered in their finite (type or dimensional)
counterparts in previous work of one of us (e.g. \cite{gln}, \cite{BeCoGaVi}, \cite{SDIWP}, 
\cite{BeGa}, \cite{GSCH}) and this paper set itself in the same stream.

%More broadly, let  $G(r,n)$ denotes the Grassmann variety 
%parametrizing $r$-dimensional subspaces of $\CC^n$. By  {\em finite dimensional}
%boson-fermion correspondence, we mean the composition of the 
%Poincar\'e duality $H^*(G(r,n);\QQ)\sra H_*(G(r,n);\QQ)$ together with the isomorphism $H_*(G(r,n);
%\QQ)\cong \bw^r\QQ^n$. The latter can be certainly understood in terms of the geometric 
%Satake correspondence, as suggested  e.g. in \cite{Labelle,ngo2021} but, more trivially, 
%noticing that their bases have the same cardinality.

\claim{}  \label{sec1:12}

To more precisely anticipate the shape of our results,  let $n\in\NN\cup\{\infty \}$, so that $V_n=\QQ[X]/(X^{n})$ is the $
\QQ$-algebra of 
all polynomials of 
degree $<n$ in one indeterminate $X$. It will be considered a $B_{1,n}:=\displaystyle{\QQ[c_1]\over (c_1^n)}$-module via the map $c_1(p(X)+(X^n))=Xp(X)+(X^n)$. Denote by $V_n^*$ its (restricted) dual, 
spanned by $(\d^j)_{0\leq j<n}$, where $\d^j$ denotes the unique linear form on $V_n$ such 
that $
\d^j(X^i)=\delta^{ij}$. If $n=\infty$ we write $V:=V_\infty=\QQ[X]$. For all $
(I,J)=(i_1,\ldots,i_h,j_1,\ldots,j_k)\in \NN^h\times\NN^k$, the formal words (free fermions)
\be
X^I\d^J:=X^{i_1}\cdots X^{i_{{h}}}\d^{j_1}\cdots\d^{j_k}\label{eq1:formwrd}
\ee 
induce homogeneous (elementary) endomorphisms of $\bw V_n$ of degree ${h-k}$ (which 
is zero if $h<k$) i.e., for all $0\leq r<n$, a  
family  of maps $\Phi_{(I,J),n}:\bw^rV_n\sra \bw^{r+h-k}V_n$  defined by
\begin{eqnarray}
\Phi_{(I,J),n}(u):&=&X^{i_1}\cdots X^{i_k}\d^{j_1}\cdots\d^{j_k}(u)\cr
&=&X^{i_1}\w\cdots\w X^{i_h}
\w \d^{j_1}\lrcorner(\ldots\lrcorner (\d^{j_h}\lrcorner u))\ldots),\quad \forall 
u\in\bw^rV_n\label{eq:wcon}
\end{eqnarray}
{ where if $\alpha\in V^*$, by $\alpha\lrcorner: \bw V_n\sra \bw V_n$ one understands the unique odd derivation of the exterior algebra such that $\alpha\lrcorner u=\alpha(u)$. It is a derivation of degrre $-1$ meaning that it maps $\bw^kV_n$ to $\bw^{k-1}V_n$ (Cf. Section~\ref{sec:sec47})}.

The exterior algebra $\bw V_n$ 
is an irreducible representation of the Clifford algebra {spanned by} the identity 
$1:=1_{\bw 
V_n}$ and all the words of type 
\eqref{eq1:formwrd}, subject to the relations $X^i\d^j+\d^jX^i=\delta^{ij}$, $X^iX^j=-
X^jX^i$ and $\d^i\d^j=-
\d^j\d^i$ holding in $\End_\QQ(\bw V_n)$.
Let $
(\mS_\blamb)_{\blamb\in\Pcal_{r,n}}$ be the Schur basis of $B_{r,n}$ (basically the Schur determinant of the Chern Classes of the universal quotient bundle over $G(r,n)$), where $\Pcal_{r,n}$ denotes the set of all partitions  
$
\blamb:=(n-r\geq 
\lambda_1\geq\cdots\geq \lambda_r\geq 0)$, in turns isomorphic to $\bw^rV_n$. Wedging with $h$ vectors and contracting with 
$k$ linear forms against elements of $\bw^rV_n$ induce obvious  vector space 
homomorphisms
\be
\Phi_{(I,J),n};B_{r,n}\lra B_{r+h-k,n}\label{eq:ntralmp}
\ee
which is our aim  to encode all together into a generating function. 

To this purpose, let 
$$
\bfz_h:=(z_1,\ldots,z_h),\qquad \bfw_k:=(w_1,\ldots, w_k)\quad  \mathrm{and} \quad \bft_r:=(t_1,
\ldots,t_r)
$$ be indeterminates over $\QQ$. For each $\blamb\in\Pcal_{r,n}$ denote by 
$s_\blamb(\bft_r)$ the Schur polynomial (Section \ref{sec:secschur}) attached to the 
partition $\blamb$ and to the indeterminates $\bft_r$, and 
$$
p_i(\bft_r):=t_1^i+\cdots+t_r^i.
$$
be the (symmetric polynomial)  sum of the $i^{th}$ powers of the indeterminates $t_j$.
Our {\bf main Theorem} \ref{thm4:mainthm} then gives an explicit 
 precise description of a $B_{r+h-k,n}$-valued formal power seriers 

\be
\Phi_{(r,h,k)}(\bfz_h,\bfw_k^{-1},\bft_r)=(-1)^h\prod_{i=1}^kw_i^{-r+1}\cdot \DD_{(r,h,k),n} 
\exp\left(\sum_{i\geq 1}x_ip_i(\bfz_h,\bft_r)+{p_i(\bft_r)\over i\prod_{j=1}^kw_j^i}\right)
\label{eq1:thmint}
\ee
which we call {\em structural generating series},
enjoying  the following property

\begin{quotation}
\noindent
The coefficient of $\boxed{\bfz_h^I\bfw_k^{-J}s_\blamb(\bft_r)}$ in the expansion of $\Phi_{(r,h,k)}
(\bfz_h,\bfw_k,
\bft_r)$ is {
the} image in $B_{r+h-k,n}$ of $\mS_\blamb\in B_{r,n}$  through the natural 
map $\Phi_{(I,J),n}$ \eqref{eq:ntralmp}. 

\end{quotation}
The key block of \eqref{eq1:thmint} is the determinant $\DD_{(r,h,k),n}$, whose explicit 
expression is displayed in \eqref{eq3:finalsh}. It has a number of relevant particular 
cases. If $h=k$ and $n=\infty,$ formula \eqref{eq1:thmint} simplifies the main statement of 
\cite{BeCoGaVi}, where one studied $B_r:=B_{r,\infty}$ as a module over the Lie algebra of 
the endomorphisms of $\bw^kV$. If in addition $h=k=1,$ one obtains  the {\em finite 
type} counterpart of the Date-Jimbo-Kashiwara-Miwa (DJKM) representation of $gl_\infty(\QQ)
$, as 
studied in \cite{gln}. Moreover,  if $r\sra \infty,$ the formula provides a determinantal 
expression of the 
vertex operators occurring in the DJKM representation (Cf. also Section \ref{sec:specialc}, Example b)). 
 The computations heavily rely on the use of the 
exponential of the 
trace representation of the Lie algebra $gl_n(\QQ)$, as in \cite{BeCoGaVi,BeGa}, which in \cite{SCHSD} was phrased in terms of {\em Hasse-Schmidt derivation} on an exterior algebra.
The general formula we obtain involves a determinant which may be considered as a deformed 
version of the usual Giambelli's one occurring in classical Schubert calculus.
Our main formula  also recovers,  as a limit when $r,n\sra \infty$, in the same spirit of previous work (e.g. 
\cite{BeCoGaVi,gln,BeCoMa}), 
the DJKM vertex operator representation of the Lie algebra $gl_\infty(\QQ)$, to emphasize that the latter is nothing but an extremal case of a 
general situation (sloganized in \cite{gln} as {\em the cohomology of the Grassmannian is a $gl_n$-module}) whose other extremal one is the usual multiplication of matrices by 
vectors.  

{The referee kindly pointed out that Theorem~\ref{thm4:mainthm} above  also admits a 
natural interpretation from the viewpoint of symmetric functions. We fully agree. Via 
the finite boson--fermion correspondence, the structural generating series acts on the Schur 
basis of $B_{r,n}$, so that the determinant $\DD_{(r,h,k)}$ may be regarded as a deformation 
of the classical Jacobi--Trudi/Giambelli determinant. This perspective is closely related to 
the philosophy developed by Jing and Rozhkovskaya \cite{JingRozh}, in which the Clifford 
algebra action is constructed starting from generalized Jacobi--Trudi identities. The present 
paper adopts a different starting point, namely the action of the Clifford algebra on the 
exterior algebra, from which the corresponding determinantal identities on the bosonic side 
are recovered via the finite boson--fermion correspondence. The two viewpoints therefore 
appear to be complementary, and we believe that it would be worthwhile to make the relationship between the two approaches more explicit in future work.}

\medskip

\noindent

% 
% In 
% case $h=k=1$ one obtains a more explicit structural generating series in the sense of 
% \cite{BeCoMa}. Our main result, however, is a very general formula, which works uniformly 
% for all $r,n$, involving Giambelli's like determinants occurring in Schubert Calculus and 
% which for $r,n\sra \infty$ tends to the description, obtained by the same methods, of the 
% Bosonic representation of the Fermionic Clifford algebra as shown in \cite{BeGa}.

%In the sequel we generalize the situation, by studying the action of the Clifford algebra $\mathrm{End}(\bw V)$ on all the rings $B_{r,n}$, which translates into a family of maps
%$$
%B_{r,n}\sra B_{r+h-k}\llb\bfz_h,\bfw_k^{-1}\rrb
%$$
% encoded in a  $B_{r+h-k}$-valued formal power series in the variables $(\bfz_h,
% \bfw_k^{-1},\bft_r)$. In case $n=\infty$, $h=1$, $k=0$ or $h=0$ and $k=1$, one obtains the 
% action of finite type bosonic vertex operators (as in e.g. \cite{gln,HSDGA}), in case 
% $h=k$ and $r<\infty$ one obtains a yet more explict descritpion of the ring $B_r$ as a 
% module over the degree $0$ endomorphisms of the exterior algebra as in \cite{BeCoGaVi}. 

\claim{\bf About the content.} The paper is organized as follows.   Section \ref{sec:sec1} 
collects a few combinatorial preliminaries which are, on the other hand, well available in 
the literature, possibly with a slight change of notation. We also discuss a couple of special cases, to show the relationship of the present subject to the theory of symmetric functions as in \cite[p.~96]{MacDonald} and/or the vertex operators occurring in the description of the KP hierarchy (generating series of Bernstein { operators like in \cite[I.5, Example~29, formula (6)]{MacDonald})}.
For the sequel, it is important to notice that 
we will be working with two rings of symmetric polynomials in two different set of variables, 
$B_r$ and $\Lambda_r$. The variables of $B_r$ are not explicitly invoked, while we use those
of the ring $\Lambda_r$. We hope this will not cause too much confusion.

Section \ref{sec:sec3} provides further details about the $B_{r,n}$-module structure of the exterior power $\bw^rV_n$. The key homomorphisms $\sigma_+(z), \ovsig_+(z):\bw V_n\sra \bw V_n\llb z\rrb$ will be introduced and  rigorously defined in terms of 
exponential of linear differential operators on the exterior algebra, and they are related with the classical Chern polynomials of the tautological and quotient bundle in the universal sequence over the Grassmannian. In this section we also show (Proposition~\ref{prop3:prop36}) how product of $r$ {\em Schubert derivations} supply a generating function for the basis elements of the $r$-th exterior power $\bw^k\QQ^n$, which is suited to study the $gl_n(\QQ)$-action of $B_{r,n}$, which is the finite dimensional shadow of the DJKM representation of $gl_\infty(\QQ)$.

In Section 
\ref{sec:sec4},  the necessary combinatorics needed to work with the representation 
of the  Clifford algebra of ``free fermions'' on an exterior algebra  is developed. This yoga was already partially 
used in \cite{BeGa, BeCoGaVi} for the bosonic representation of the Clifford algebra 
generalizing formulas due \DJKM \,\,\cite{DJKM01}. 

Section \ref{sec:sec5} is finally devoted to 
the statement and proof of our main formula. As anticipated, it involves a determinant that is as huge as the 
integer $r+h$, but whose expression has a very elegant and regular shape. Using Mathematica$^\copyright$,  Example \ref{ex6:ex65} fully works out the expression of the generating structural series of the action of the $4\times 4$ square matrices against the cohomology  of the Grassmannian $G(2,4)$: it is the first instance of the DJKM representation in a finite dimensional context which is not just the multiplication of a matrix times a vector.
\paragraph{\bf Acknowledgment.}  We are grateful to Andr\'e Contiero, Frank Neumann, Renato Vidal Martins, Parham Saleyhan  for useful comments and thoughtful remarks, and Edoardo Ballico for constant encouragement. The present version of the paper has greatly benefited from the careful reading and valuable suggestions of the anonymous referee, whose comments substantially improved its exposition.

The first author was partially supported by the PRIN Multilinear Algebraic Geometry, the italian INDAM--GNSAGA, and Ricerca di Base of Politecnico di Torino (grant no. 53$\_$RBA23GATLET). He also
thanks Politecnico di Torino for covering the open-access publication of the present article.

The second author thanks ICTP (Trieste) for supporting her visit to Politecnico di Torino, where part of this work was carried out in a stimulating and welcoming research environment. She is also grateful to Professor Rasoul Roknizadeh for helpful discussions on the physical background of the subject. This research was partially supported by IPM grant No. 1403170414, was partly carried out at the IPM–Isfahan Branch, and was also supported by the ICTP Research in Pairs Programme (Trieste) and by the Iran National Science Foundation (INSF, grant No. 4001480).

\paragraph{Data Availability Statement.} The authors declare that the data supporting the findings of this study are available within the paper. All the article listed in the final reference list are publicly available either on journals or in the {\tt ArXiv} data repository.

\section{Preliminaries and {Notation}}\label{sec:sec1}
\claim{} If $A$ is any algebra over some commutative ring $R$ and $\bft_r:=(t_1,\ldots,t_r)$ are indeterminates over $A$, the notation $A\llb \bft_r\rrb$ stands for the algebra of formal power series in those indeterminates with the usual definition of product.
\claim{} Let $V:=\QQ[X]$ be the vector space of polynomials with rational coefficients. Its natural basis is $\bfX:=(X^i)_{i\geq 0}$.  For each $n\in \NN\cup\{\infty\}$,  let 
\begin{center}$V_n:=\displaystyle{\QQ[X]\over(X^{n-1})}$ \quad if $n<\infty$\quad  and \quad $V_\infty:=V$.
\end{center} The linear maps 
$$
\d^i:={1\over i!}\left.{d^i\over dX^i}\right|_{X=0}:V_n\mapsto \QQ
$$
define the linear forms ${\bm\d}:=(\d^i)_{i\geq 0}$ of the (restricted) dual $V_n^*$ of $V_n$, dual to $\bfX=(X^i)_{0\leq i<n}$, i.e. $\d^j(X^i)=\delta^j_i$. 
We denote by $\bfX(z)$ and $\bd(w^{-1})$ the generating functions of the bases $\bfX$ and $\bd$ of $V$ and $V^*$ respectively, i.e.
\be
\bfX(z)=\sum_{0\leq i<n}X^iz^i\qquad \mathrm{and}
\qquad \bd(w^{-1})=
\sum_{0\leq i< n}\d^iw^{-i}=\exp\left(w^{-1}{d\over dX}\right).\label{eq2:genfnd}
\ee

\claim{} Let $(r,n)\in (\NN\cup\{\infty\})^2$ such that $0\leq r\leq n$. To any such pair 
we attach three isomorphic vector spaces, as follows. The first one is  $\PP_{r,n}:
\QQ^{\Pcal_{r,n}}$, the $\QQ$-vector space freely generated by all the {\em partitions} 
of length at most $r$ whose Young diagram is contained in a $r\times (n-r)$ rectangle; 
the second one is $B_{r,n}$, the universal decomposition algebra of the polynomial $T^n$ 
as the product of two monic polynomials of degree $r$ and $n-r$ (see e.g. \cite{BeNa}, 
\cite{LakTh01,LakTh02}), and $\bw^rV_n$, the $r$-th exterior power of an $n$-dimensional 
vector space. In all of the above, $n$ may be infinite. The three spaces are related by 
the commutative diagram
\be
\begin{tikzcd}[row sep=2em, column sep=2em]
& \mathbb{P}_{r,n} \arrow[dl, "\beta"'] \arrow[dr, "\phi"] & \\
B_{r,n} \arrow[rr, "\phi\circ\beta^{-1}"', below] & & \bigwedge\nolimits^r V_n
\end{tikzcd}\label{eq:bfc_cd}
\ee
where $\phi\circ \beta^{-1}$ can be understood as the {\em boson-fermion} correspondence, as will be explained in the sequel. 

\claim{\bf The space $\PP_{r,n}$.} A {\em partition} is a {non-increasing} sequence $\blamb:
\lambda_1\geq\lambda_2\geq\cdots\geq  0$ of non negative integers with finite length $
\ell(\blamb):=\sharp\{i\,|\, \lambda_i\neq 0\}$. We denote by $\Pcal$ the set of all 
partitions, by $\Pcal_r:=\{\blamb\in\Pcal\,|\,\ell(\blamb)\leq r\}$ the set of all partitions 
of length at most $r$ and by $\Pcal_{r,n}:=\{\blamb\in\Pcal_r\;|\, \lambda_1\geq n-r\}$ the 
set of all partitions whose Young diagram is contained in a $r\times (n-r)$ rectangle.  
Clearly $\Pcal_r=\Pcal_{r,\infty}$ and $\Pcal=\Pcal_\infty$. Let $
\PP_{r,n}:=\QQ^{\Pcal_{r,n}}$ be the $\QQ$-vector space generated by all 
$\blamb\in \Pcal_{r,n}$
 Each element  $u\in \PP_{r,n}$ can be written in {\em bosonic}  ($\beta(u)$) or {\em fermionic} ($\phi(u)$) notation, in the sense specified below.
 
\claim{\bf Exterior powers and algebra.} \label{sec:secschur} 
The $r$-th exterior power of $V_n$ is the vector space $\bw^rV_n$ generated by all the expressions of the form $X^{i_1}\w\cdots\w X^{i_r}$, modulo the relations 
$$
X^{i_\tau(1)}\w\cdots \w X^{i_\tau(n)}=\sgn(\tau)X^{i_1}\w\cdots\w X^{i_r},\quad \tau\in S_n
$$ where $\sgn(\tau)$ denotes the sign of the permutation $\tau$.  Therefore, 
$
\bw^rV_n=\bigoplus_{\blamb\in\Pcal_{r,n}}\QQ\cdot\bfX^r(\blamb)
$,
where
\be
\bfX^r(\blamb)=\bfX^r(\lambda_1,\ldots,\lambda_r):=X^{r-1+\lambda_1}\w\cdots\w X^{\lambda_r},
\ee
whence the $\QQ$-linear isomorphism $\phi:\PP_{r,n}\sra \bw^rV_n$ mapping $\blamb$ to $\bfX^r(\blamb)$.

The {\em exterior algebra} of $V_n$ is $(\bw V_n,\w)$,  where $\bw V_n=\bigoplus_{r\geq 0}\bw^rV_n$ and $$\w: \bw^rV_n\times \bw^sV_n\sra \bw^{r+s}V_n$$ is the (bilinear) the juxtaposition product.

	 Consider  now the wedge product
$$
\bfX(t_1)\w\cdots\w \bfX(t_r)\in (\bw^rV_n)\llb t_1,\ldots,t_r\rrb.
$$
It  is clearly a (possibly infinite) linear combination of $(\bfX^r(\blamb))_{\blamb\in\Pcal_{r,n}}$ which is 
divisible by the Vandermonde determinant $\Delta_0(\bft_r):=\prod_{i<j}(t_i-t_j)$. The 
equality
$$
\bfX(t_1)\w\cdots\w \bfX(t_r)=\Delta_0(\bft_r)\sum_{\blamb\in\Pcal_{r,n}}
\bfX^r(\blamb)s_\blamb(\bft_r)
$$
defines the {\em Schur symmetric polynomials} {$s_\blamb(\bft_r)$} in $\bft_r=(t_1,\ldots,t_r)
$. Therefore
$$
\sum_{\blamb\in\Pcal_{r,n}}\bfX^r(\blamb)s_\blamb(\bft_r)\in{ \big(\bw^rV_n\big)}\llb t_1,\ldots, t_r\rrb
$$
is a generating function of the basis elements of $\bw^rV_n$ \footnotemark\footnotetext{{Here $(\bw V_n)\llb t_1,\ldots, t_r \rrb$ are the $\bw V_n$-valued formal power series in the 
indeterminates $(t_1,\ldots, t_r)$}.}. It is well known that  $\big(s_\blamb(\bft_r)
\big)_{\blamb\in\Pcal_{r}}$ forms  a basis of the ring $\Lambda_r:=\QQ[\bft_r]^{S_r}$ of the 
symmetric polynomials in the indeterminates $(t_1,\ldots, t_r)$.
 
\claim{\bf The ``bosonic space'' $B_{r,n}$.}\label{sec:sec21} Let $\bfc:=(c_1,c_2,\ldots)$ be 
a sequence of indeterminates over $\QQ$. To each  $r\in\NN\cup\{\infty\}$ we attach a $
\QQ$-algebra $B_r$.  By definition $B_0=\QQ$ and, for 
$0<r\leq \infty$, one set  $B_r:=\QQ[\mc_1,\ldots, \mc_r]$. In addition one sets $B:=B_\infty=\QQ[\bfmc]$ for the polynomial $\QQ$-
algebra in the infinitely many indeterminates $(\mc_i)$. One considers the generic 
(Chern) monic polynomial\footnotemark\footnotetext{The indeterminate coefficients are 
denoted $c_i$, to reminish the Chern classes of the tautological bundle over the 
Grassmannian $G(r,n)$.}
\vspace{-20pt}
$$
\bfmc_r(z)=1-\mc_1z^r+\cdots+(-1)^r\mc_rz^r\in B_r[z]
$$
and the formal (Segre) power series \footnotemark\footnotetext{The indeterminate coefficients are denoted $S_i$, to reminish the Segre classes of the tautological bundle over the Grassmannian $G(r,n)$.}
$$
\bfmS_r(z)=\sum_{{i}\in \ZZ}\mS_iz^i:={1\over \bfmc_r(z)}\in  B_r\llb z\rrb.
$$

In other words :
$$
\bfmc_r(z)(1+\mS_1z+\cdots+\mS_{n-r}z^r)=1
$$
the equality holding in $B_{r,n}[z]$ for finite $n$, or
$$
\bfmc_r(z)(1+\sum_{j\geq 1}\mS_jz^j)=1
$$
the equality holding in $B_r\llb z\rrb$  if $n=\infty$.

For $n\geq r\geq 0$, let $I_{r,n}:=(\mS_{n-r+1},\ldots, \mS_{n})$ be  the ideal of $B_{r}$ generated by $\mS_{n-r+1},\ldots, \mS_{n}$. Let $$
\pi_{r,n}:B_r\lra B_{r,n}={B_r\over I_{r,n}}. 
$$
be the canonical projection.
It is well known that $B_{r,n}$ is the universal decomposition algebra of the monomial $T^n$ as the product of two monic polynomials, one of degree $r$ (see e.g. \cite{BeNa}). This is in turn the presentation of the singular cohomology ring 
of the Grassmannian, once the indeterminates $(\mc_i)$ are interpreted as the Chern classes of the
tautological rank $r$ bundle  and the {$\mS_i$'s} as the Chern classes
of the universal quotient bundle over the Grassmannian $G(r,n)$ (or the Seegre classes of the tautological bundle).
 We denote $\bfmS_{r,n}=(\pi_{r,n}(\mS_j))_{j\in\ZZ}=(1=\mS_0,\mS_1,\ldots, \mS_{n-r})$ and
$$
\bfmS_{r,n}(z)=1+\mS_1z+\cdots+\mS_{n-r}z^r\in B_{r,n}[z]
$$
%which through the identification of $B_{r,n}$ with the cohomology $H^+(G(r,n),\QQ)$ would correspond to the {\em Segre polynomial.}
Let $\bfx_{r,n}=(x_1,x_2,\ldots)$  be the sequence of elements of $B_{r,n}$ defined via the equality
\be
\exp\left(\sum_{i\geq 1}x_iz^i\right)=\bfmS_{r,n}(z)=1+S_1z+\cdots+\mS_{n-r}z^{n-r}.\label{eq22:xi}
\ee
In other words
$$
\bfx_{r,n}(z):=\sum_{i\geq 1}x_iz^i=\log(1+\mS_1z+\cdots+\mS_{n-r}z^{n-r})=\sum_{j\geq 1}{1\over j}\Big(S_1z+\cdots+\mS_{n-r}z^{n-r}\Big)^{j-1}.
$$

Let $p_i(\bft_r):=\sum_{j=1}^r t_j^i$. To each partition $\blamb\in \Pcal_{r,n}$ we associate $S_\blamb\in\Pcal_{r,n}$ defined as follows
\be
\sum_{\blamb\in\Pcal_r}\mS_\blamb\cdot  s_\blamb(\bft_r)=\exp(\sum_{i\geq 1}x_ip_i(\bft_r)).\label{eq3:preJJTT}
\ee
It is well known that $(S_\blamb\,|\,\blamb\in\Pcal_{r,n})$ is a $\QQ$-basis of $B_{r,n}$:
$$
B_{r,n}:=\bigoplus_{\blamb\in\Pcal_{r,n}}\QQ\cdot \mS_\blamb
$$
so that $\exp\left(\sum_{i\geq 1}x_ip_i(\bft_r)\right)$ is a generating function of the 
basis $S_\blamb$. 
This gives us the second alluded natural isomorphisms $\beta:\PP_{r,n}\sra B_{r,n}$, mapping $\blamb$ to $S_\blamb$.

Moreover, the equality
\be
\exp\left(\sum_{i\geq 1}x_ip_i(\bft_r)\right)=\exp\left(\sum_{\lambda_1\geq 1}
x_it_1^{\lambda_1}\right)\red{\cdots}\exp\left(\sum_{\lambda_r\geq 1}
x_i\bft_r^{\lambda_r}\right),\label{eq2:pJT}
\ee
 along with \eqref{eq3:preJJTT}, implies the celebrated Jacobi--Trudy formula
$$
\mS_\blamb=\det(\mS_{\lambda_j-j+i})
$$
obtained by equating the coefficients of $s_\blamb(\bft_r)$ on both members of \eqref{eq2:pJT}.
Since the $\mS_i$'s generate $B_{r,n}$ as a ring, the structure of $B_{r,n}$ as a $\QQ$-algebra is determined once one knows all the products
$$
\mS_i\cdot\mS_\blamb,
$$
as linear combination of the basis $(\mS_\blamb)$. All such products are ruled by the obvious equality
\be
\exp\left(\sum_{i\geq 1}x_iz^i\right)\exp\left(\sum_{i\geq 1}x_ip_i(\bft_r)\right)=\exp\left(\sum_{i\geq 1}x_ip_i(z,\bft_r)\right)=\sum_{\blamb\in\Pcal_{r,n}}\mS_\blamb\cdot s_\blamb(z,\bft_r).\label{eq2:genpieri}
\ee
In other words
$
\mS_i\cdot S_\blamb 
$
is the coefficient of $z^is_\blamb(\bft_r)$ in the most right hand side of \eqref{eq2:genpieri} or, more practically, the coefficient of $z^i\bft_r^{\blamb}$ 
of
$$
\Delta_0(\bft_r)\cdot \exp\left(\sum_{i\geq 1}x_ip_i(z,\bft_r)\right),
$$ 
where $\bft_r^\blamb:=t_1^{r-1+\lambda_1}t_2^{r-2+\lambda_2}\cdots t_r^{\lambda_r}$.
The map $\beta: \PP_{r,n}\sra B_{r,n}$ maps  any $\blamb \in \PP_{r,n}$ to its bosonic expression $S_\blamb$.

\claim{} 
%Our mains actors are the cohomology rings $B_{r,n}$, with 
%rational coefficients,  of the complex Grassmannians $G(r,n)$, parametrizing 
%$r$-dimensional vector subspaces of $\CC^n$. We simply write $B_r$ if $n=\infty$, which is the cohomology of the ind-manifold $G(r,\infty)$. 
If both $r=n=\infty$, 
 one is then concerned with the polynomial ring $B:=\QQ[x_1,x_2,\ldots,]$ in infinitely many 
indeterminates. It corresponds to the cohomology ring of the Universal Grassmann Manifold (UGM), an {\em ind-variety}  which Sato constructs in \cite{sato1} as a direct limit of finite dimensional Grassmannians. As the UGM is a direct limit of 
finite dimensional Grassmannians, its cohomology is the inverse limit of the 
cohomologies of the finite dimensional ones, and coincides with $B_\infty$.

\section{Generating functions of linear maps: a few special cases}\label{sec:specialc}

Recall that the purpose of this paper is to provide, inspired by the  DJKM vertex operator representation of the Lie algebra $gl(\infty)$, explicit descriptions of maps $B_{r,n}\sra B_{r+s.n}$ induced by elementary endomorphisms of the exterior algebra $\bw V_n$, where $V_n=\QQ[X]/(X^n)$.

It seems convenient  to anticipate a few special cases which, though elementary, seems interesting in their own. 
\begin{enumerate}
\item[a)]
 If $r=1$, $B_{1,n}
=\displaystyle{\QQ[\mc_1]\over{(\mc_1^n)}}$.  It  is obviously isomorphic to $V_n=\QQ[X]/(X^{n})$. We think of $V_n$ as a $B_{1,n}$-module such that $c_1X^i=X^{i+1}+(X^{n})$. 

Any 
matrix $A\in \mathrm{End}_\QQ(V_n)$ can be seen as linear combination of $E_{ij}:=X^i\otimes \d^j$ (notation as in \ref{sec1:12}). Let  $\Ecal_{1,n}=\sum_{i,j=0}^nX^iz^i\otimes\d^jw^{-j}$ which can be seen as the $\QQ[z,w^{-1}]$ valued matrix
$$
\begin{pmatrix} 1& w^{-1}&\cdots& w^{-n+1}\cr
z& zw^{-1}&\cdots& zw^{-n+1}\cr
\vdots&\vdots&\ddots&\vdots\cr
z^{n-1}&z^{n-1}w^{-1}&\cdots&z^{n-1}w^{-n+1}
\end{pmatrix}
$$
If 
$X^k$ is any basis element of $V_n$ one has
$$
\Ecal_{1,n}(z,w)X^k=\sum_{i=1}^kX^iz^i\delta^j_k=\sum_{i=1}^k\mc_1^iX^0z^iw^{-k}={w^{-k}\over 1-\mc_1z}X^0
$$  
The expression can be further improved if, instead of acting just on $X^k$, we act on all the 
$X^k$ at once ($1\leq k\leq n$), using a generating series, which provides the elegant formula below:
\be
\Ecal_{1,n}(z,w,t)X^0:=\Ecal_{1,n}(z,w)\sum X^kt^k={1\over 1 -\mc_1z}\left({1-\displaystyle{t^n\over w^n}
\over 1-\displaystyle{t\over w}}\right)X^0={1-\displaystyle{t^n\over w^n\displaystyle}\over (1-\mc_1z)
\left(1-\displaystyle{t\over  w}\right)}\,X^0\label{eq1:byhand}
\ee
Since working over $\QQ$, we have the exponential form:
\begin{eqnarray}
\Ecal_{1,n}(z,w,t)&=&\left(1-\displaystyle{t^n\over w^n}\right)\exp\left(\sum_{i\geq 1}{1\over i}
\left(\mc_1^nz^n+\displaystyle{t^i\over w^i}\right)\right)\cr
&=&\exp\left(\sum_{i\geq 1}{1\over i}
\left(\mc_1^nz^n+\displaystyle{t^i\over w^i}-\displaystyle{t^{ni}\over w^{ni}}\right)\right).\label{eq1:ssfps}
\end{eqnarray}
We call \eqref{eq1:ssfps} the {\em generating structural series} governing the multiplication of a 
matrix by a vector of $V_n$, identified with a polynomial of degree $\leq n-1$. In other words
$$
E_{ij}X^k=\Big([z^iw^{-j}t^k]\Ecal_{1,n}(z,w,t)\Big)X^0.
$$
Expression \eqref{eq1:ssfps}  is  a particular case of the general formula \eqref{eq4:formnthm} for $r=h=k=1$ and $n<\infty$. If $n=\infty$ one has
$$
\Ecal_{1}(z,w,t):=\Ecal_{1,\infty}(z,w,t)=\exp\left(\sum_{i\geq 1}{1\over i}
\left(\mc_1^nz^n+\displaystyle{t^i\over w^i}\right)\right),
$$
another case taken into account by the same formula \eqref{eq4:formnthm}.

\item[b)] The second example we propose points to the relationship of our subject  to vertex operators and  the classical theory of symmetric functions at once. For sake of simplicity, let $n=\infty$, so working in the ring $B_r=\QQ[\mc_1,\ldots,\mc_r]
$.  The task is to describe the action of wedging a vector against a polynomial of 
$B_r$. To do so, one recalls that the basis element $\mS_\blamb$ of $B_r$ ($
\blamb\in\Pcal_r$) does correspond to the basis element $X^{r-1+\lambda_1}\w\cdots\w X^{\lambda_r}  \in
\bw^rV$. We are so interested in computing $X^i\w \mS_\blamb$, defined by
\be
(X^i\w \mS_\blamb)\bfX^{r+1}(0):=X^i\w (\mS_\blamb\bfX^{r}(0)):= X^i\w \bfX^r(\blamb)= X^i\w X^{r-1+\lambda_1}\w\cdots\w X^{\lambda_r}.\label{eq1:wedgbos}
\ee
%Instead of wedging $\mS_\blamb$  just by $X^i$, we do by wedging against the generating series $\sum_{i\geq 0}
%X^iz^i$. 
Let us now stipulate the existence, established in Section \ref{sec3:secsig}, of an 
algebra homomorphism
$$
\sigma_+(z):\bw V\sra \bw V\llb z\rrb
$$
such that 
$$
\sigma_+(z)X^0=\sum_{i\geq 0}X^iz^i={X^0\over 1-Xz}.
$$
Denote by $\ovsig_+(z)$ the inverse of $\sigma_+(z)$ in the power series algebra $\End_\QQ(\bw  V)\llb z\rrb$: it is a $\w$-homomorphism as well, i.e. $\ovsig_+(z)(u\w v)=\ovsig_+(z)u\w\ovsig_+(z)v$, and is the unique such that $\ovsig_+(z)u=u-(Xu)z$, for all $u\in V$. Therefore
\begin{eqnarray}
\hskip-10pt\sigma_+(z)X^0\w \big(X^{r-1+\lambda_1}\w\cdots\w X^{\lambda_r}\big)&=&\sigma_+(z)\Big(X^0\w \ovsig_+(z)(X^{r-1+\lambda_1}\w \cdots\w X^{\lambda_r})\Big)\cr\cr
&=& \sigma_+(z)(X^0\w \ovsig_+(z)X^{r-1+\lambda_1}\w\cdots\w \ovsig_+(z) X^{\lambda_r} \label{eq1:for}
\end{eqnarray}
Formula \eqref{eq1:for} can be recast in two different ways. The former consists in considering the unique $\w$-algebra homomorphism $\ovsig_-(z):\bw V\sra \bw V\llb z^{-1}\rrb$ (like in \cite{gln,SDIWP}) such that
$$
\ovsig_-(z)_{|V}=1-X^{-1}z^{-1},
$$
where by $X^{-1}$ one means the (locally) nilpotent endomorphism of $V_n$ defined by $X^{-1}X^j=X^{j-1}$ if $j\geq 1$ and $0$ otherwise. Then 
 \eqref{eq1:for} can be written as
\beqns
\sigma_+(z)X^0\w X^{r-1+\lambda_1}\w\cdots\w X^{\lambda_r}&=&z^r\sigma_+(z)\ovsig_{-}(z)(X^{r+\lambda_1}\w\cdots  X^{1+\lambda_r}\w X^0)\cr
&=&z^r\Gamma_r(z)S_{(\lambda_1,\cdots,\lambda_r)}\bfX^{r+1}(0)
\eqns
where one sets
$$
(\Gamma_r(z)\mS_{(\lambda_1,\cdots,\lambda_r)})X^{r+1}\w\cdots\w X^1\w X^0:=\sigma_+(z)\ovsig_{-}(z)\bfX^{r+1}(\lambda_1,\cdots,\lambda_r).
$$
It is not that difficult to see, as shown in \cite{HSDGA,pluckercone}, that
$$
\Gamma(z):=\lim_{r\sra \infty} \Gamma_r(z)=\exp\left(\sum_{i\geq 1}x_iz^i\right)\exp\left(\sum_{i\geq 1}{1\over iz^i}{\d\over d x_i}\right)
$$
which is precisely the expression of one of the  Bernstein (vertex) operators occurring in the description of the Boson-Fermion correspondence (\cite{KacRaRoz,MJDbook} or \cite[p.~96]{MacDonald}). See also Rosas \cite{Rosas} for further combinatorial aspects of these vertex operators.

This paper takes however an alternative approach, which brings  the subject back in the 
realm of the classical theory of symmetric functions, where vertex  operators are clues of the
presence of  Vandermonde determinants. 
%Indeed, the above exercise gets simpler once we evaluate the wedging action not just against one basis element of $\bw^rV_n$ but on the generating series of all of them.
Let $\sigma_+(\bft_r)=\sigma_+(t_1)\cdots\sigma_+(t_r)$ (a kind of multivariate Hasse-Schmidt derivation on the
exterior algebra, see also \cite{fereshteh}). One of our outputs is the 
remarkable equality (Proposition \ref{prop3:prop36})
$$\sum_{\blamb\in\Pcal_2}\bfX^r(\blamb)s_\blamb(\bft_r)=\sigma_+(\bft_r)\left(X^{r-1}
\w\cdots\w X^0\right)
$$ 

This time one wedges $\sigma_+(z)X^0$  against not just one basis element of $\bw^rV$, but 
rather the generating series of the basis of $\bw^rV$. 

%where $s_\blamb(\bft_2)$ is the Schur symmetric polynomial living in the ring $\Lambda_2$ of symmetric polynomials in the indeterminates $\bft_2=(t_1,t_2)$,  is achieved by the formula
%$$
%\sigma_+(\bft_2)X^1\w X^0=\sigma_+(t_1)\sigma_+(t_2)X^1\w X^0=\sigma(t_1,t_2)X^1\w \sigma_+(t_1,t_2)X^0.
%$$
In this case we have that since 
$$
\sigma_+(z)X^0\w \sigma_+(\bft_r)\left[X^{r-1}\w\cdots\w X^0\right]
$$
vanishes whenever $z=t_i$ for some $1\leq i\leq r$, it can be written as
$$
\prod_{i=1}^r(z-t_i)\sigma_+(z,t_1,\cdots,t_r)\left[X^r\w X^{r-1}\w\cdots\w X^1\w X^0\right].
$$
From a  ``bosonic'' point of view, i.e. with $B_r$ coefficients,  it can be seen as 
\be
z^r\prod_{i=1}^r\left(1-{t_i\over z}\right)\exp\left(\sum_{i\geq 1}x_ip_i(z,\bft_r)
\right)=z^r\exp\left(\sum_{i\geq 1}{p_i(\bft_r)\over z^i}\right)\exp\left(\sum_{i\geq 1}
x_ip_i(z,\bft_r)\right),
\ee
%
%z^2\exp\left(\sum_{i\geq 1}x_i(z^i+t_1^i+t_2^i)\right)X^2\w X^1\w X^0=
%$$
which is one more particular case of the main formula \eqref{eq4:formnthm}, putting $h=1$ and $k=0$.
In passing, we also remark  that
\begin{eqnarray*}
z^r\exp(\sum_{i\geq 1}x_iz^i)\Gamma_r(z)\exp\left(\sum_{i\geq 1}{x_ip_i(\bft_r)}\right)
%&=&z^r\exp\left(\sum_{i\geq 1}{x_ip_i(\bft_r)\over z^i}\right)\exp\left(\sum_{i\geq 1}x_ip_i(z,\bft_r)\right)\cr\cr
=z^r\exp(\sum_{i\geq 1}x_iz^i)\exp\left(\sum_{i\geq 1}\left(x_i+{1\over z^i}\right)p_i(\bft_r)\right)
\end{eqnarray*}
whence 
\be
\Gamma_r(z)\exp\left(\sum_{i\geq 1}{x_ip_i(\bft_r)}\right)=\exp\left(\sum_{i\geq 1}\left(x_i+{1\over z^i}\right)p_i(\bft_r)\right).\label{eq1:berns}
\ee
For $r=\infty$, formula~\eqref{eq1:berns} is nothing but the sum over all the partitions of the right hand side of the formula stated in  \cite[Theorem 3.4]{CarGou}. It is a consequence of  the  commutation rule (similar to those proven within a similar framework see \cite[Section 4]{BeCoGaVi} or \cite[Section 8]{SDIWP}).
$$
\ovsig_-(z)\sigma_+(\bft_r)=\left(\sum_{i\geq 1}{p_i(\bft_r)\over z^i}\right)\sigma_+(\bft_r)\ovsig_-(z),
$$
occurring also in a more general categorical context like in \cite[formula (12), (13),  (14)]{FrPeSe}.
In sum,  wedging a generating function of the basis of $\bw^hV$  against a generating function of the basis of $B_r$, in the sense of \eqref{eq1:wedgbos}, is the same as multiplying two exponentials. 
\end{enumerate}

\section{The $B_{r,n}$--module structure of the exterior algebra}\label{sec:sec3}
To generalize the examples of Section~\ref{sec:specialc}, we  recall the natural module structure of the exterior algebra of a countably infinite dimensional $\QQ$--vector space over the ring $B_{r,n}$.

%\claim{} For each $\blamb\in \Pcal_{r,n}$, define $\mS_\blamb$ through the equality
%\be
%\sum_{\blamb\in\Pcal_{r,n}}\mS_\blamb\cdot s_\blamb(\bft_r)=\exp\left(\sum_{i\geq 1}
%x_ip_i(\bft_r)\right)\label{eq3:preJJTT}.
%\ee

\claim{}\label{sec3:secsig} Let $gl_n(\QQ):=(\QQ^{n\times n},{[\cdot,\cdot]})$ be the Lie algebra of all the $\QQ$-valued $n\times n$ square matrices with respect to the usual commutator.  We understand it as
\be
gl_n(\QQ):=\bigoplus_{0\leq i,j<n}\QQ\cdot E_{ij}\label{eq2:glnq}
\ee
where, as costumarily, we denote by $E_{ij}$ the elementary  $n\times n$ matrix whose entries are all zero but $1$ in position $(i,j)$. Therefore, due to \eqref{eq2:glnq}, if $n=\infty$, by $gl_\infty(\QQ)$ we understand the Lie algebra of {all the matrices} of infinite size whose entries are all zero but finitely many.

The exterior algebra $\bw V_n=\bigoplus_{r\geq 0}\bw^rV_n$ is a $gl_n(\QQ)$-module via the {\em trace representation}.  The  trace $\tr(A)$  of  $A\in gl_n(\QQ)$ is the unique derivation of the exterior algebra such that it restrict to the multiplication by $A$ to the first degree of $\bw V_n$, i.e.:
$$
\left\{\begin{matrix}
\tr(A)(u\w v)&=&\tr(A)u\w v+u\w\tr(A)v\cr
\tr(A)u&=&Au,\quad \forall u\in V_n=\bw^1V_n.
\end{matrix}\right.
$$

For all $j\geq 0$, let $X^j:V_n\sra V_n$ denote the endomorphism given by multiplication by $X^j$, which for finite $r$ can be identified with  the elementary $n\times n$ nilpotent matrix of rank $n-j$, and let us denote by $\delta(X^j)$ its trace.
For  any indeterminate $z$, let:
$$
\sigma_+(z)=\sum_{i\geq 0} \sigma_iz^i:=\exp\left(\sum_{j\geq 1}{\delta(X^j)\over j}z^j\right): \bw V_n\sra \left(\bw V_n\right)\llb z\rrb.
$$
Notice that $\bfX(z)=\sigma_+(z)X^0$. Because it is the exponential of a derivation, it is clearly a $\bw$-algebra homomorphism, namely $\sigma_+(z)(u\w v)=\sigma_+(z)u\w \sigma_+(z)v$.
We also use the notation 
$$
\ovsig_+(z)=\sum_{i\geq 0}(-1)^i\ovsig_iz^i=\exp\left(-\sum_{j\geq 1}{\delta(X^j)\over j}z^j\right): \bw V_n\sra \left(\bw V_n\right)\llb z\rrb
$$
which is clearly the inverse of $\sigma_+(z)$ in the algebra $\End_\QQ(\bw V_n)\llb z\rrb$, i.e.
$$
\sigma_+(z)\ovsig_+(z)=\ovsig_+(z)\sigma_+(z)=1,
$$
implying the following Cayley-Hamilton relations (see e.g. \cite{GSCH}):
\begin{eqnarray}
\sigma_+(z)(\ovsig_+(z)u\w v)&=&u\w\sigma_+(z)v,\\ \cr
\ovsig_+(z)(\sigma_+(z)u\w v)&=&u\w\ovsig_+(z)v.
\end{eqnarray}
We also denote
$$
\sigma_+(\bft_r)=\prod_{i=1}^r\sigma_+(t_i)=\sigma_+(t_1)\cdots\sigma_+(t_r),
$$
the product in  $\End_\QQ(\bw V_n)\llb \bft_r\rrb$ regarded as an algebra over $\End_\QQ(\bw V_n)$.  We have:
\bclm{\bf Proposition.}\label{prop3:prop36} {\em For all $r\geq 1$, the following equality holds:
\be  
\sigma_+(t_1)X^0\w\cdots\w \sigma_+(t_r)X^0=\Delta_0(\bft_r)\sigma_+(\bft_r)\bfX^r(0).\label{eq3:fakVa}
\ee
}
\eclm
\proof For $r=1$, we set $\Delta_0(t_1)=1$, so that formula \eqref{eq3:fakVa} holds in this case. For $r=2$ we have, using the Cayley Hamilton relations :
\begin{eqnarray*}
\sigma_+(t_1)X^0\w \sigma_+(t_2)X^0&=&\sigma_+(t_1,t_2)\big(\ovsig_+(t_2)X^0\w\ovsig_+(t_1)X^0\big)\label{eq214:ch24}\cr\cr
&=&\sigma_+(t_1,t_2)\big[(X^0-t_1X^1)\w (X^0-t_2X^1)\big]\cr\cr
&=&(t_1-t_2)\sigma_+(t_1,t_2)X^1\w X^0
\end{eqnarray*}
and thus the property holds for $r=2$. Then we argue by induction. Let us suppose that 
$$
\sigma_+(t_1)X^0\w\cdots\w\sigma_+(t_{r-1})X^0=\Delta_0(\bft_{r-1})\sigma_+(\bft_r)\bfX^{r-1}(0).
$$
Then, on one hand $\sigma_+(t_1)X^0\w\cdots\w\sigma_+(t_{r-1})X^0\w\sigma_+(t_r)X^0$ is clearly a multiple of $\prod_{i<j}(t_i-t_j)$. On the other hand, using induction:
\begin{eqnarray*}
\sigma_+(t_1)X^0\w\cdots\w\sigma_+(t_{r-1})X^0\w\sigma_+(t_r)X^0=\Delta_0(\bft_{r-1})\sigma_+(\bft_{r-1})\bfX^{r-1}\w \sigma_+(t_r)X^0
\end{eqnarray*}
Invoking the Cayley-Hamilton relations \eqref{eq214:ch24} we get, from the last side:
\begin{eqnarray}
&&\Delta_0(\bft_{r-1})\sigma_+(\bft_{r-1})\bfX^{r-1}\w \sigma_+(t_r)X^0\cr\cr
&=&\Delta:0(\bft_{r-1})\sigma_+(\bft_r)\ovsig_+(t_r)\bfX^{r-1}(0)\w \ovsig_+(\bft_r\hskip-2pt\setminus t_r)X^0\cr\cr
&=&\sum_{i=0}^{r-1}(-1)^it^i\bfX^{r-1}(1^i)\w (X^0-e_1(\bft_{r-1}X^1+\cdots+(-1)^{r-1}e_{r-1}(\bft_{r-1})X^{r-1}.\cr&&\label{eq3:redvndm}
\end{eqnarray}
The above expression is clearly a multiple of $\bfX^r(0)$, which is divisible by $\prod_{i=1}^r (t_i-t_r)$. To determine its coefficient one looks at that of $t_1\cdots t_{r-1}=e_{r-1}(\bft_{r-1})$ in \eqref{eq3:redvndm}. In all cases one gets a summand of the form:
\be
(-1)^{r-1}t_1\cdots t_{r-1} \bfX^{r-1}(0)\w X^{r-1}.\label{eq3:exch}
\ee
Now:  
$$
\bfX^{r-1}(0)\w X^{r-1}=(-1)^{r-1} X^{r-1}\w\bfX^{r-1}(0)=(-1)^{r-1}\bfX^r(0).
$$
 Therefore
$$
(-1)^{r-1}e_{r-1}(\bft_{r-1})=(-1)^{r-1}t_1\cdots t_{r-1} \bfX^{r-1}(0)\w X^{r-1}=t_1\cdots t_{r-1}\bfX^r(0).
$$

proving that \eqref{eq3:redvndm} is equal to $\Delta_0(\bft_r)\bfX^{r}(0)$, as claimed.\qed
%\bclm{\bf Proposition.}\label{prop3:prop36} {\em
%$$
%\sigma_+(t_1)X^0\w\cdots\w \sigma_+(t_r)X^0=\Delta_0(\bft_r)\cdot\sigma_+(\bft_r)\bfX^r(0)=\Delta_0(\bft_r)\sum_{\blamb\in\Pcal_{r,n}}\bfX^r(\blamb)s_\blamb(\bft_r)
%$$
%where $\Delta_0(\bft_r)=\prod_{1\leq i<j\leq r}(t_i-t_j)$  is the Vandermonde determinant. Therefore $\sigma_+(\bft_r)\bfX^r(0)$ is the generating function of the basis $(\bfX^r(\blamb))_{\blamb\in\Pcal_{r,n}}$ of $\bw^rV_n$.}
%\eclm
%\proof
%Let $\bft_r\setminus t_i:=\{t_1,\ldots,\widehat t_i,\ldots,t_r\}$ ($t_i$ omitted). Then
%$$
%\sigma_+(t_1)X^0\w\cdots\w \sigma_+(t_r)X^0=\sigma_+(\bft_r)\left(\ovsig_+(\bft_r\setminus t_1)X^0\w\cdots\w\ovsig_+(\bft_r\setminus t_r)X^0\right)
%$$
%which is clearly  $P(\bft_r)\bfX^r(0)$, where $P(\bft_r)$ is a polynomial of degree  $(r-1)(r-2)\cdots 1=(r-1)!$, vanishing whenever $t_i=t_j$ it is divisible by $(t_i-t_j)$ for all $r\geq j>i\geq 1$. Thus $P(t)$  is a constant multiple of the Vandermonde $\Delta_0(\bft_r)=\prod_{i<j}(t_j-t_i)$. To determine the constant is sufficient to look at the coefficient of $t_r^{r-1}t_{r-1}^{r-2}\cdots t_2$ in the expansion of
%$$
%\ovsig_+(\bft_r)X^0\w\cdots\w\ovsig_+(\bft_r)X^0
%$$
%which is $(-1)^{r-1}$.

\bclm{\bf Proposition.} \label{pro37:prop37} {\em The equality
\[
x_i\cdot \bfX^r(\blamb):=\delta(X^i)\bfX^r(\blamb),\qquad \blamb\in\Pcal_{r,n}
\]
 makes $\bw^r V_n$ into a free $B_{r,n}$-module of rank $1$ generated by $\bfX^r(0)=X^{r-1}\w X^{r-2}\w\cdots\w X^0$, such that $\bfX^r(\blamb)=\mS_\blamb\cdot\bfX^r(0)$.} 
\eclm
\proof Since $\mS_i$ is a polynomial in $\bfx_{r,n}=(x_1,x_2,\ldots)$, it is enough to show
that
$$
\mS_i(\mS_ju)=(\mS_i\mS_j)u,\qquad {(\forall u\in \bw^rV_n)}.
$$
We will check it through generating functions:
\beqns
\exp\left(\sum_{i\geq 1}x_iz^i\right)\exp\left(\sum_{i\geq 1}x_iw^i\right)\bfX^r(\blamb)&=&\exp\left(\sum_{i\geq 1}x_iz^i\right)\exp\left(\sum_{j\geq 1}{\delta(X^j)\over j}w^j\right)\bfX^r(\blamb)\cr\cr =\exp\left(\sum_{i\geq 1}{\delta(X^i)\over i}z^i\right)\exp\left(\sum_{j\geq 1}{\delta(X^j)\over j}w^j\right)\bfX^r(\blamb)
&=&\exp\left(\sum_{i\geq 1}{\delta(X^i)\over i}(z^i+w^i)\right)\bfX^r(0)\cr
=\exp\left(\sum_{i\geq 1}x_i(z^i+w^i)\right)\bfX^r(0)
&=& \left(\exp(\sum_{i\geq 1}x_iz^i)\exp(\sum_{i\geq 1}x_iw^i)\right)\bfX^r(0).
\eqns
The assignment $x_i\mapsto {\delta(X^i)\over i}\bfX^r(0)$ is clearly injective. To prove that it is surjective, it is enough to show that $\bfX^r(\blamb)=\mS_\blamb\cdot \bfX^r(0)$ and we do that all at once by using generationg functions:
\beqns
\sum_{\blamb\in\Pcal_{r,n}}\bfX^r(\blamb)s_\blamb(\bft_r)&=&\sigma_+(\bft_r)\bfX^r(0)=\exp\left(\sum_{i\geq 1}{\delta(X^i)\over i}p_i(\bft_r)\right)\bfX^r(0)\cr\cr
&=&\exp\left(\sum_{i\geq 1}x_ip_i(\bft_r)\right)\bfX^r(0)=\sum_{\blamb\in\Pcal_{r,n}}\mS_\blamb\bfX^r(0)s_\blamb(\bft_r).
\eqns
More explicitly, we have the equality:
\be
\bfX^r(\blamb)=X^{r-1+\lambda_1}\w\cdots\w X^{\lambda_r}=\begin{vmatrix}
\mS_{\lambda_1}&\mS_{\lambda_2-1}&\cdots&\mS_{\lambda_r-r+1}\cr
\mS_{\lambda_1+1}&\mS_{\lambda_2}&\cdots&\mS_{\lambda_r-r+2}\cr
\vdots&\vdots&\ddots&\vdots\cr
\mS_{\lambda_1+r-1}&\mS_{\lambda_2+r-1}&\cdots&\mS_{\lambda_r}
\end{vmatrix}\cdot X^{r-1}\w X^{r-2}\w\cdots\w X^0.\label{eq2:JTfext}
\ee
which can also be seen as the action of the word $X^{r-1+\lambda_1}\w\cdots \w X^{\lambda_r}$ on $B_0=\QQ$!
\bclm{\bf Corollary.} {\em The isomorphism $B_{r,n}\mapsto \bw^rV_n$ makes $B_{r,n}$ into a representation of the Lie algebra $gl_n(\QQ)$.}
\eclm
%The isomorphism  $P\mapsto P\cdot X^r(0)$ defined in Proposition 3.7 together with the  trace action of $gl_n(\QQ)$ on $\bw^rV_n$, makes $B_{r,n}$ into a representation of the Lie algebra $gl_n(\QQ)$.

\proof Indeed, the action $\circledast$ defined by  \be
(A\circledast P)\bfX^r(0)=\tr(A)(P\cdot\bfX^r(0)).
\ee
satisfies $[A,B]\circledast P=A\circledast(B\circledast P)-B\circledast(A\circledast P)$ because $\tr([A,B])=[\tr(A),\tr(B)]$, as well known.\qed
\claim{} We recall the formalism of two-rows-determinant introduced in \cite{BeCoGaVi} to deal with contractions. Let $(a_0,\ldots, a_{r-1})$ be any set of indeterminates (or rational numbers) and $(v_0, v_1,\ldots, v_{r-1})\in V_n^r$. Then we let
\be
\begin{vmatrix}
a_0&a_1&\cdots&a_{r-1}\cr\cr
v_0&v_1&\ldots&v_{r-1}
\end{vmatrix}=\sum_{j=0}^{r-1}{(-1)^ja_j\bfv_{I\setminus\{j\}}}\label{eq2:2rwsdet}
\ee
where  {$I=\{0,\ldots,r-1\}$ and $\bfv_{I\setminus\{j\}}$ is the wedge product of $v_0,v_1,\ldots, v_{r-1}$ with $v_j$ removed.}
For example

$$
\begin{vmatrix}
a_0&a_1&a_2\cr\cr
v_0&v_1&v_2
\end{vmatrix}=a_0\,v_1\w v_2-a_1\,v_0\w v_2+a_2\,v_0\w v_1.
$$

\smallskip
\noindent
\bclm{\bf Remark.}\label{rmk:rmklin} Notice that the  the output of expression  \eqref{eq2:2rwsdet} changes sign whenever two columns are exchanged. Hence, the linearity with respect to the first entry of the second row implies the  linearity with respect to all the others arguments.
\eclm

\claim{} \label{sec:sec47} The {\em contraction}  of a vector of $\bw^rV_n$ against a linear form $\alpha$ can be { defined} using diagram \eqref{eq2:2rwsdet} { by declaring}:
$$
\alpha\lrcorner v:=\begin{vmatrix}
\alpha(v)\cr v
\end{vmatrix}=\alpha(v)
$$
and, for $r\geq 2$
\be
\alpha\lrcorner (v_1\w v_2\w\cdots\w v_r)=\begin{vmatrix}
\alpha(v_1)&\cdots&\alpha(v_r)\cr\cr
v_1&\cdots&v_r
\end{vmatrix}.\label{eq2:contdt}
\ee
\bclm{\bf Proposition.}\label{prop2:gmbcnt}  {\em Let $r\leq n-1$. For each $(a_1,\ldots,a_r)$ {and a partition $\lambda_1\geq \cdots\geq \lambda_r$},  the following formula holds
\be
\begin{vmatrix}
a_1&a_2&\cdots&a_r\cr\cr
X^{r-1+\lambda_1}&X^{r-2+\lambda_2}&\cdots&X^{\lambda_r}
\end{vmatrix}=\begin{vmatrix}
a_1&a_2&\cdots&a_r\cr \mS_{\lambda_1+1}&\mS_{\lambda_2}&\cdots&\mS_{\lambda_r+r-2}\cr
\vdots&\vdots&\ddots&\vdots\cr
\mS_{\lambda_1+r-1}&\mS_{\lambda_2+r-2}&\cdots&\mS_{\lambda_r}
\end{vmatrix}X^{r-2}\w\cdots\w X^0\label{eq2:redJT}
\ee
}
\eclm
\proof It is a consequence of \eqref{eq2:JTfext}. Given a partition $\blamb=(\lambda_1\geq\cdots\geq\lambda_r)$, let $\blamb^{(j)}$ be the partition of length at most $r-1$ obtained {by} removing $\lambda_j$ from $\blamb$. Then
\beqns
\begin{vmatrix}
a_1&a_2&\cdots&a_r\cr\cr
X^{r-1+\lambda_1}&X^{r-2+\lambda_2}&\cdots&X^{\lambda_r}
\end{vmatrix}&=&\sum_{i=1}^r(-1)^{j-1}a_j\bfX^{r-1}(\blamb^{(j)})\cr\cr\cr
&=&(\sum_{i=1}^r(-1)^{j-1}a_j\mS_{\blamb^{(j)}})\bfX^{r-1}(0)
\eqns
which is precisely the expansion of the determinant occurring on the right hand side of \eqref{eq2:redJT}.\qed
\bclm{\bf Corollary.} \label{cor2:mncor}{\em  Let $\bfw_k:=(w_1,\ldots,w_k)$ be a $k$-tuple of indeterminates and  {keep the same} notation as in \eqref{eq2:genfnd}. Consider the map
$$
\d(\bfw_k^{-1})\lrcorner :\bw^rV_n\sra \bw^{r-k}V_{n}\llb w_1^{-1}, \ldots,w_k^{-1}\rrb
$$
defined by 
$$
\d(\bfw_{k}^{-1})\lrcorner u=\d(\bfw_{k-1}^{-1})\lrcorner(\d(w_1^{-1})\lrcorner u).
$$  
%$$
%u\mapsto \lrcorner\bfX^r(\blamb):=\bd(\bfw_k^{-1})\lrcorner\bfX^r(\blamb):=\bd(\bfw_k^{-1})\lrcorner\left(\bd(\bfw_{k-1}^{-1})
%\lrcorner \cdots \lrcorner \d(\bfw_1^{-1}\right)\lrcorner u
%$$
Then

$$
\bd(\bfw_k^{-1})\lrcorner\bfX^r(\blamb):=\bd(w_k^{-1})\lrcorner\left(\bd(w_{k-1}^{-1})
\lrcorner \cdots \lrcorner \d(w_1^{-1})\lrcorner \bfX^r(\blamb)\right)
$$

\bigskip
$$
 =\begin{vmatrix}
w_1^{-r+1-\lambda_1}&w_1^{-r+2-\lambda_2}&\cdots&w_1^{-\lambda_r}\cr
\vdots&\vdots&\ddots&\vdots\cr
w_k^{-r+1-\lambda_1}&w_k^{-r+2-\lambda_2}&\cdots&w_k^{-\lambda_r}\cr\cr
X^{r-1+\lambda_1}&X^{r-2+\lambda_2}&\cdots&X^{\lambda_r}
\end{vmatrix}
$$

\bigskip
$$
=w_1^{-r+1}w_2^{-r+2}\cdots w_k^{-r+k}\begin{vmatrix}
w_1^{-\lambda_1}&w_1^{1-\lambda_2}&w_1^{2-\lambda_3}&\cdots&w_1^{r-1+\lambda_r}\cr
\vdots&\vdots&\vdots&\ddots&\vdots\cr
w_k^{1-k-\lambda_1}&w_k^{2-k-\lambda_2}&w_k^{-k+3-\lambda_3}&\cdots&w_k^{r-k+\lambda_r}\cr\cr
\mS_{\lambda_1+k}&\mS_{\lambda_2+k-1}&\mS_{\lambda_3+k-2}&\cdots&\mS_{\lambda_r+k+r-3}\cr
\vdots&\vdots&\ddots&\vdots\cr
\mS_{\lambda_1+r-1}&\mS_{\lambda_2+r-2}&\mS_{\lambda_2-r+1}
&\cdots&\mS_{\lambda_r}
\end{vmatrix}X^{r-1-k}\w\cdots\w X^0.
$$
}
\eclm

\smallskip
\proof
One first notice that
$
\bd(w_i^{-1})X^j=w_i^{-j}
$
and then  merges it in the rephrasing \eqref{eq2:contdt} of the contraction, 
eventually using Proposition~\ref{prop2:gmbcnt}. \qed
\section{The action of the Clifford Algebra}\label{sec:sec4}
Let $V_n^*:={\bigoplus_{j= 0}^{n-1}}\QQ\cdot\d^j$ be the restricted dual of $V_n$ 
(which, for $n<\infty$ is the usual dual). Let us consider the canonical 
Clifford Algebra $\Ccal(V)$ supported on $V\oplus V^*$, defined by the canonical 
non-degenerate bilinear form
$$
\langle u\oplus\alpha,v\oplus\beta \rangle=\beta(u)+\alpha(v).
$$
Write $\bd^{r}(\bmu)=\d^{r-1+\mu_1}\w \cdots\d^{\mu_r}\in \bw^rV_n^*$.
An element (a word) $\bfX^h(\blamb)\bd^k(\bmu)$ of the Clifford algebra $
\Ccal(V)$, in {\em normal form}, acts on the exterior algebra $\bw V$ by 
wedging and contracting:
\be
\bfX^h(\blamb)\bd^k(\bmu)(u){:}=\bfX^h(\blamb)\w (\bd^k(\bmu)\lrcorner u)\qquad (\forall u\in\bw V).
\label{eq:general}
\ee
If an element is not in normal form, one uses the well known and easy to check relations:
$$
\left[X^i,\d^j\right]_+\red{:}=X^i\d^j+\d^jX^i=\delta_{ij}{\mathbf 1}_{\bw V}.
$$
The $X^i$ and $\d^j$ plays the same role as the generators of the free-fermions Clifford algebra, usually denoted by $\psi_i, \psi_i^*$ (see e.g. \cite{jimbomiwa} or \cite[Section 1]{SavageBF}).
{Notice that if $n<\infty$, the Clifford algebra is naturally isomorphic to the 
algebra of the endomorphisms $gl(\bw V_n)$ of $\bw V_n$, isomorphic to $\bw 
V_n\otimes \bw V_n^*$.} If $n=\infty$ the action is different, because the image 
of an element of the form $u\otimes 1$ is in general an infinite linear combination of the basis elements of $\bw V_n$, contrarily to the image 
of any element of $u\otimes \alpha$, with $\alpha\neq 1$.
\claim{} Let $\bfz_h=(z_1,\ldots,z_h)$ and $\bfw_k^{-1}=(w_1^{-1},\ldots, 
w_k^{-1})$ be, respectively, an $h$-tuple and a {$k$-tuple of indeterminates}. 
Let
$$
\Ecal_{(r,h,k);n}(\bfz_h,\bfw_k^{-1}):=\sum_{(\blamb,\bmu)\in \Pcal_{h,n}\times 
\Pcal_{k,n}}\bfX^h(\blamb)\otimes 
\bd^k(\bmu)s_\blamb(\bfz_h)s_\bmu(\bfw_k^{-1}).
$$
Thus, for each $r,h,k\geq 0$, the Clifford algebra $\Ccal(V)$ induces 
vectorspace homomorphisms
$$
B_{r}\lra B_{r+h-k}\llb\bfz_h,\bfw_k^{-1}\rrb
$$
defined by
$$
(\Ecal_{h,k,n}(\bfz_h,\bfw_k^{-1}))\bfX^{r}(0)=\sum_{\blamb,\bmu}
\bfX^h(\blamb)\w \big(\bd^k(\bmu)\lrcorner \bfX^r(0)
\big)s_\blamb(\bfz_k)s_\bmu(\bfw_k^{-1})
$$
which is clearly zero if $r+h-k<0$.
This enables us to define a graded linear homomorphims 
\be 
\star:B_{r,n}\sra B_{r+h-k,n}\label{eq:starr}
\ee
 through the equality:
$$
\bfX^h(\blamb)\bd^k(\bmu)\star P:=\left[s_\blamb(\bfz_h)s_\bmu(\bfw_k^{-1})
\right]\left(\Ecal_{(r,h,k),n}(\bfz_h,\bfw_k^{-1})P\right)\in \bw^{r+h-k}
V_n\llb\bfz^k,\bfw_k^{-1}\rrb.
$$

Recall that $E_{i,j}\in gl_n(\QQ)$ denotes the elementary matrix with all entries zero but 
$1$ in position $(i,j)$, acting on $u\in V_n$ as
$
E_{i,j}(u)=(X^i\otimes\d^j)(u)=X^i\d^j(u)
$
\bclm{\bf Proposition.}\label{prop3:tra} For all $u\in \bw V_n$,
$
\tr(E_{ij})u= X^i\w (\d^j\lrcorner u)
$
\eclm
\proof
If $u\in V_n$, then $\tr(E_{ij})(u)=E_{ij}(u)=X^i\w \d^j(u)$ and the property holds. Let us assume that it holds for each $v\in \bw^jV_n$ for all $0\leq j\leq r-1$. Then if $u\in V_n$ and $v\in \bw V_n$ is arbitrarily chosen:
\beqns
\tr(E_{ij})(u\w v)&=&\tr(E_{ij})u\w v+u\w \tr(E_{ij})v=X^i\d^j(u)\w v+ u\w X^i\w (\d^j\lrcorner v)\cr
&=&X^i\d^j(u)\w v-X^i\w \d^j\lrcorner v=X^i\w (\d^j\lrcorner (u\w v)).
\eqns
Because of each element of $\bw V_n$ is a finite linear combination of monomials of the form $u\w v$, with $u\in V_n$, the property is proven.
\qed

If $h=k=1$, put  $\bfz_1=z$ and $\bfw_1^{-1}=w^{-1}$. Then $\Ecal_{r,n}(z,w):=\Ecal_{(r,1,1);n}(z,w)$ is  the generating function of the  elementary endomorphisms of $V_n$:
$$
\Ecal_{r,n}(z,w):=\sum_{0\leq i,j<n}E_{i,j}z^iw^{-j}\in gl_n(\QQ)[[z,w^{-1}]].
$$

$$
\Ecal_{r,n}(z,w)=\sum_{0\leq i,j<n}X^iz^iw^{-j}\otimes \d^jw^{-j}=\bfX(z)\otimes \bd(w^{-1})=\sigma_+(z)X^0\otimes \bd(w^{-1})
$$
where $\bd(w^{-1})=\sum_{i\geq 0}\d^iw^{-i}$ is a generating function of the dual basis of $\bfX$. 
It follows that the action of $\Ecal_{r,n}(z,w)$ is given by
\[
\Big(\Ecal_n(z,w)\star S_\blamb)\Big)\bfX^r(0)=\tr(\Ecal_n(z,w))\bfX^r(\blamb)
\]
a particular case of \eqref{eq:general}.

\bclm{\bf Proposition.} For all $u\in \bw V_n$
$$
\tr(\Ecal_{r,n}(z,w))u=\sigma_+(z)X^0\w \bd(w^{-1})\lrcorner u
$$
\eclm
\proof
{It is a} straightforward consequence of Proposition \ref{prop3:tra}.\qed 

\section{The Main Theorem}\label{sec:sec5}
\claim{} Our main goal is to explicitly compute the action of $\Ecal_{(r,h,k);n}(\bfz_h,\bfw_k^{-1})$ against the generating function $\exp(\bfx_{r,n}(\bft_r))$ of the $\QQ$-basis of $B_{r,n}$.
Recall that (Proposition \ref{prop3:prop36})
$$
\exp(\bfx_{r,n}(\bft_r))\bfX^r(0)=\sigma_+(\bft_r)\bfX^r(0),
$$
To this purpose, 
%
%Our goal is now that of computing the action of $\Ecal_{r,n}(z,w)$ against the generating function of a basis of $\bw^rV_n$.
%\bclm{\bf Proposition}
%$$
%\sum_{\blamb\in\Pcal_{r,n}}X^r(\blamb)s_\blamb(\bft_r)=\sigma_+(\bft_r)\bfX^r(0)=\exp\left(\sum_{i\geq 1}x_ip_i(\bft_r)\right)\bfX^r(0)
%$$
%\eclm
%\proof
%It follows from the definition of the module structure. \qed
one needs one more piece of notation.

\claim{} Let ${\mathfrak S}_r$ be the symmetric group on $r$-elements and $\bft_r=(t_1,
\ldots,t_r)$ be an $r$-tuple of {indeterminates} and let
$$\Lambda_r:=\QQ[\bft_r]^{{\mathfrak S}_r}$$  be the ring  of symmetric polynomials in $\bft_r$. 
For an arbitrary {indeterminate} $\zeta$ over $\QQ[\bft_r]$, let 
$$
E_r(\bft_r,\zeta){:=}1-e_1(\bft_r)\zeta+\cdots+(-1)^re_r(\bft_r)\zeta^r=\prod_{i=1}^r(1-
t_i\zeta)\in \QQ[\bft_r],
$$
which is the generating function of the {\em elementary symmetric polynomials}. The equality
$$
H_r(\bft_r,\zeta)=\sum_{j\in\ZZ}h_j(\bft_r)\zeta^j={1\over E_r(\bft_r,\zeta)}
$$
defines the complete symmetric polynomials $h_j(\bft_r)$ in $\bft_r$. Let
$$
p_i(\bft_r):=t_1^i+\cdots+t_r^i
$$
be the $i$-th symmetric Newton power sum. They are related to $e_i(\bft_r)$ and $h_j(\bft_r)
$ via the equality 
$$
\exp\left(\sum_{i\geq 1}{p_i(\bft_r)\over i}\zeta^i\right)={1\over E_r(\bft_r,\zeta)}
=\sum_{j\in\ZZ}h_j(\bft_r)\zeta^j.
$$
holding in the algebra $\Lambda_r\llb\xi \rrb$.

By the fundamental theorem of symmetric functions, $\Lambda_r{=}\QQ[e_1(\bft_r),
\ldots, e_r(\bft_r)]$. The ring $\Lambda_r$ is also generated as a $\QQ$-algebra by $
(h_1(\bft_r),\ldots,h_r(\bft_r))$, so that $h_j(\bft_r)$ is a weighted homogeneous 
polynomial of degree $j$ in $h_1(\bft_r), \ldots, h_r(\bft_r)$. The same holds for $(p_1(\bft_r),\ldots, 
p_r(\bft_r))$, i.e. they are algebraically independent and each $p_i(\bft_r)$ is a 
homogeneous polynomial expression of weighted degree $i$ in $(p_j(\bft_r))_{1\leq j\leq r}$.
\bclm{\bf Definition.} {\em For all $(j,m)\in\NN\times \NN$, define symmetric polynomials $U_{j,m}(\bft_r)$ via the equality:
$$
\sum_{j\geq 0}U_{j,m}(\bft_r)\zeta^j=E_r(\bft_r,\zeta)\sum_{i\geq m}h_i(\bft_r)\zeta^i.
$$
}
\eclm
 To get a feeling of how they explicitly 
look like is easy.
For instance
$$
U_{0,m}(h_m(\bft_r))=h_m(\bft_r),\qquad U_{1,m}(h_m(\bft_r))=h_{m+1}(\bft_r)-e_1(\bft_r)h_m(\bft_r),
$$
$$
 U_{2,m}(\bft_r)=h_{m+2}(\bft_r)-e_1(\bft_r)h_{m+1}(\bft_r)+e_2(\bft_r)h_m(\bft_r),
$$
In general: \quad
$
U_{j,m}(\bft_r)=h_{m+j}-e_1(\bft_r)h_{m+j-1}(\bft_r)+\cdots+(_1)^re_r(\bft(r)h_m(\bft_r)
$

\bclm{\bf Proposition}\label{prop1:finsum}{\em
\begin{enumerate}
\item For all $m\geq 0$, $U_{r+j,m}(\bft_r)=0$ for all $j\geq 0$;
\item The sum of the first $m$ terms of $H_r(\bft_r,\zeta)$ is
$$
\sum_{j=0}^{m-1}h_j(\bft_r)\zeta^j={1\over E_r(\bft_r,\zeta)}(1-y_m\zeta^m)
$$
where
\be
y_m=U_{0,m}(\bft_r)+U_{1,m}(\bft_r)\zeta+\cdots+U_{r-1,m}(\bft_r)\zeta^{r-1}\label{eq1:shtcut}
\ee
\end{enumerate}
}
\eclm

\proof
To prove the first part, notice that for all $j\geq 0$
$$
U_{r+j,m}(\bft_r)=h_{m+r+i}(\bft_r)-e_1(\bft_r)h_{m+r+i-1}(\bft_r)+\cdots+(-1)^rh_{m+i}(\bft_r)
$$
which vanishes due to the relation $E_r(\bft_r,\zeta)H_r(\bft_r,\zeta)=1$ (all the coefficients of the positive powers of the product vanish). The second part is a consequence of the first one. We have
\beqns
1&=&E_r(\bft_r,\zeta)\sum_{j=0}^{m-1}h_j(\bft_r)\zeta^j+E_r(\bft_r,\zeta)\sum_{j\geq m}h_j(\bft_r)\zeta^j\cr\cr\cr
&=&E_r(\bft_r,\zeta)\sum_{j=0}^{m-1}h_j(\bft_r)\zeta^j+\zeta^m\Big(U_{0,m}(\bft_r)+U_{1,m}(\bft_r)\zeta+\cdots+U_{r-1,m}(\bft_r)\zeta^{r-1}\Big)
\eqns
due to the vanishing of $U_{j,m}(\bft_r)$ for all $j\geq r$. Therefore
$$
\sum_{j=0}^{m-1}h_j(\bft_r)\zeta^j
%={1\over E_r(\bft_r,\zeta)}(1-\zeta^m(U_{0,m}(\bft_r)+U_{1,m}(\bft_r)\zeta+\cdots+U_{r-1,m}%(\bft_r)\zeta^{r-1})
={1\over E_r(\bft_r,\zeta)}(1-y_m\zeta^m)=\exp\left(\sum_{j\geq 1}{p_j(\bft_r)\over j}\zeta^j\right)(1-y_m\zeta^m)
$$
where $y_m$ is as in \eqref{eq1:shtcut}. \qed

\claim{} One still needs to define  linear combinations $\mScal_j(\bft_r)\in B_{r,n}\otimes\Lambda_r$ of the polynomials 
$
\mS_j\in B_{r,n}
$
with coefficients in $\Lambda_r$. They are  defined through the equality
\be
\sum_{j\geq -r}\mScal_j(\bft_r)\zeta^{j}:=E_r(\bft_r,\zeta^{-1})\sum_{j\geq 0}\mS_j\zeta^{j}.\label{eq0:bothm}
\ee
Equating the coefficients of the same power of {$\zeta$} on both members of \eqref{eq0:bothm}, one obtains the first values of $\mScal_i(\bft_r)$, and the feeling of what they look like:
$$
\mScal_{-r}(\bft_r)=1;\qquad \mScal_{-r+1}(\bft_r)=-e_1(\bft_r)\mS_0+e_2(\bft_r)\mS_1;\quad \mScal_0(\bft_r)=\mS_0-e_1(\bft_r)\mS_1+\cdots+(-1)^re_r(\bft_r)\mS_r
$$ 
or
\[
{\mScal}_n(z)=\mS_n-z\mS_{n+1}, \qquad
\mScal_n(z,w)=\mS_n-(z+w)\mS_{n+1}+zw\cdot \mS_{n+2}.
\]
Again,  keeping into account that $S_j=0$ for $j<0$:
\[
\mScal_{-2}(\bft_r)=e_2(\bft_r)-e_3(\bft_r)\mS_1+\cdots+(-1)^{r+1}e_r(\bft_r)\mS_{r-2}.
\]
\claim{} Let us finally consider the $(r+h)\times (r+h)$ square matrix given by:
$$
\hskip-2pt\DD_{(r,h,k),n}\hskip-2pt=\hskip-2pt
$$
\begin{small}
\be
\begin{vmatrix}
0&\cdots&0&1-\displaystyle{y_{n-r+1}w_1^{-n+r-1}}&w_1\left(1-\displaystyle{y_{n-r+2}w_1^{-n+r-2}}\right)&\hskip-8pt \cdots&\hskip-10pt w_1^{r-1}\left(1-\displaystyle{y_{n}w_1^{-n}}\right)\cr\cr
\vdots&\ddots&\vdots&\vdots&\vdots&\ddots&\vdots\cr\cr
0&\cdots&0&1-\displaystyle{y_{n-r+1}w_k^{-n+r-1}}&w_k\left(1-\displaystyle{y_{n-r+2}w_k^{-n+r-2}}\right)&\hskip-8pt \cdots&\hskip-10pt w_k^{r-1}\left(1-\displaystyle{y_{n}w_k^{-n}}\right)\cr\cr\cr
\mScal_{-r+k-1}(\bft_r)&\cdots&\mScal_{-r+k-h}(\bft_r)&\mScal_{-k+1}(\bfz_h)&\mScal_{-k}(\bfz_h)&\hskip-8pt\cdots&\mScal_{-r+k}(\bfz_h)\cr\cr
\mScal_{-r+k-2}(\bft_r)&\cdots&\mScal_{-r+k-h+1}(\bft_r)&\mScal_{-k+2}(\bfz_h)&\mScal_{-k+1}(\bfz_h)&\hskip-8pt\cdots&\mScal_{-r+k+1}(\bfz_h)\cr\cr
\vdots&\ddots&\vdots&\hskip-8pt\vdots&\vdots&\ddots&\vdots\cr\cr
\mScal_{h-1}(\bft_r)&\cdots&\mScal_0(\bft_r)&\mScal_{r-1}(\bfz_h)&\mScal_{r-2}(\bfz_h)&\cdots&\mScal_0(\bfz_h)
\end{vmatrix}\label{eq3:finalsh}
\ee
\end{small}

Recalling that (Cf. \eqref{eq:starr} for the definition of $\star$)

$$
\left(\Ecal_{(r,h,k),n}(\bfz_h,\bfw_k)\star\exp(\sum_{i\geq 1}x_ip_i(\bft_r))\right)\bfX^r(0):=
\Ecal_{(r,h,k),n}(\bfz_h,\bfw_k)\sigma_+(\bft_r)\bfX^r(0)
$$
we are now in position to state our main theorem.

\bclm{\bf Theorem.}\label{thm4:mainthm} 
$$
\Ecal_{(r,h,k),n}(\bfz_h,\bfw_k)\star\exp(\sum_{i\geq 1}x_ip_i(\bft_r))
$$
\be
=(-1)^h\prod_{i=1}^kw_i^{-r+1}\cdot \DD_{(r,h,k),n} \exp\left(\sum_{i\geq 1}x_ip_i(\bfz_h)+{p_i(\bft_r)\over i\prod_{j=1}^kw_j^i}\right)\exp(\sum_{i\geq 1}x_ip_i(\bft_r))\label{eq4:formnthm}
\ee
\eclm
\proof
The first step of the proof consists in using the equality
\beqns
\Ecal_{(r,h,k),n}(\bfz_h,\bfw_k)&\star&\exp(\sum_{i\geq 1}x_ip_i(\bft_r))\bfX^r(0)=\sigma_+(\bfz_h)\bfX^r(0)\w \d(\bfw_k^{-1})\lrcorner \sigma_+(\bft_r)\bfX^r(0)\cr\cr
&=&\sigma_+(\bfz_h)\bfX^r(0)\w \begin{vmatrix}
\d(w_1^{-1})(\sigma_+(\bft_r)X^{r-1})&\cdots &\d(w_1^{-1})(\sigma_+(\bft_r)X^0)\cr\vdots&\ddots&\vdots \cr
\d(w_k^{-1})(\sigma_+(\bft_r)X^{r-1})&\cdots &\d(w_k^{-1})(\sigma_+(\bft_r)X^0)\cr\cr

\sigma_+(\bft_r)X^{r-1}&\cdots&\sigma_+(\bft_r)X^{0}
\end{vmatrix}
\eqns
from which
\beqns
\Ecal_{(r,h,k),n}(\bfz_h,\bfw_k)&\star&\exp(\sum_{i\geq 1}x_ip_i(\bft_r))\bfX^r(0)\cr\cr\cr
&=&(-1)^h\begin{vmatrix}
0&\cdots&0&\d(w_1^{-1})(\sigma_+(\bft_r)X^{r-1})&\cdots&\d(w_1^{-1})(\sigma_+(\bft_r)X^{0})\cr
\vdots&\ddots&\vdots&\vdots&\ddots&\vdots\cr
0&\cdots&0&\d(w_k^{-1})(\sigma_+(\bft_r)X^{r-1})&\cdots&\d(w_k^{-1})(\sigma_+(\bft_r)X^{0})\cr\cr
\sigma_+(\bfz_h)X^{h-1}&\cdots&\sigma_+(\bfz_h)X^{0}& \sigma_+(\bft_r)X^{r-1}&\cdots&\sigma_+(\bft_r)X^{0}
\end{vmatrix}
\eqns

\medskip
\be
=(-1)^h\sigma_+(\bfz_k)\sigma_+(\bft_r)\begin{vmatrix}
0&\cdots&0&\d(w_1^{-1})(\sigma_+(\bft_r)X^{r-1})&\cdots&\d(w_1^{-1})(\sigma_+(\bft_r)X^{0})\cr
\vdots&\ddots&\vdots&\vdots&\ddots&\vdots\cr
0&\cdots&0&\d(w_k^{-1})(\sigma_+(\bft_r)X^{r-1})&\cdots&\d(w_k^{-1})(\sigma_+(\bft_r)X^{0})\cr\cr
\ovsig_+(\bft_r)X^{h-1}&\cdots&\ovsig_+(\bft_r)X^{0}& \ovsig_+(\bfz_h)X^{r-1}&\cdots&\ovsig_+(\bfz_h)X^{0}
\end{vmatrix}.\label{eq3:pregen1}
\ee
Now, for all $0\leq j\leq r-1$ and $1\leq i\leq k$
\beqns
\bd(w_i^{-1})(\sigma_+(\bft_r)X^{r-1-j})&=&\bd(w_i^{-1})(\sum_{0\leq i<n-r+j+1}h_i(\bft_r)X^{r-1-j+i})\cr\cr\cr
&=&\sum_{0\leq i<n-r+j+1}h_i(\bft_r)w_i^{-r+1+j-i}\cr\cr\cr
&=&w_i^{-r+1+j}\sum_{0\leq i<n-r+j+1}h_i(\bft_r)w_i^{-i}\cr\cr\cr
&=&w_i^{-r+1+j}{1\over E_r(\bft_r,w_i^{-1})}\left(1-y_{n-r+j+1}w_i^{-n+r-j-1}\right)
\eqns
Keeping into account that
$$
{1\over E_r(\bft_r,w_i^{-1})}=\exp(-\log(E_r(\bft_r,w_i^{-1}))=\exp\left(\sum_{j\geq 1}{p_j(\bft_r)\over j w^j_i}\right)
$$

and that 
$$
\sigma_+(\bft_r)\bw^rV_n=\exp\left(\sum x_ip_i(\bft_r)\right)\bw^rV_n.
$$ {Formula}  \eqref{eq3:pregen1} can be written as
$$
(-1)^h\prod_{i=1}^kw_i^{-r+1}\exp\left(\sum_{j\geq 1}x_jp_j(\bfz_h)+{p_j(\bft_r)\over j\prod_{i=1}^rw_i^j}\right)\exp\left(\sum_{i\geq 1}x_jp_j(\bft_r)\right)\widetilde{D}_{(r,h,k),n}
$$
where  using Proposition~\ref{prop1:finsum}.
\begin{large}
$$
\widetilde{D}_{(r,h,k),n}=
$$
\end{large}
\be
\begin{vmatrix}
0&\cdots&0&1-\displaystyle{y_{n-r+1}\over w_1^{n-r+1}}&w_1\left(1-\displaystyle{y_{n-r+2}\over w_1^{n-r+2}}\right)&\cdots&w_1^{r-1}\left(1-\displaystyle{y_{n}\over w_1^{n}}\right)\cr\cr
\vdots&\ddots&\vdots&\vdots&\vdots&\ddots&\vdots\cr\cr
0&\cdots&0&1-\displaystyle{y_{n-r+1}\over w_k^{n-r+1}}&w_k\left(1-\displaystyle{y_{n-r+2}\over w_k^{n-r+2}}\right)&\cdots&w_k^{r-1}\left(1-\displaystyle{y_{n}\over w_k^{n}}\right)\cr\cr
\ovsig_+(\bft_r)X^{h-1}&\cdots&\ovsig_+(\bft_r)X^{0}& \ovsig_+(\bfz_h)X^{r-1}&\ovsig_+(\bfz_h)X^{r-2}&\cdots&\ovsig_+(\bfz_h)X^{0}
\end{vmatrix}\label{eq2:prefr}
\ee

\noindent

We are only left to express $\widetilde{D}_{(r,h,k),n}$ in the shape proclaimed in the statement. To this purpose we first notice that
$$
\ovsig_+(\bft_r)X^{h-j}=X^0-e_1(\bft_r)X^1+\cdots+(-1)^re_r(\bft_r)X^{h-j+r}
$$
and, for all $0\leq j\leq r-1$
$$
\ovsig_+(\bfz_h)X^{r-j}=X^{r-j}-e_1(\bfz_h)X^{r-j+1}+\cdots+(-1)^he_h(\bfz_h)X^{r-j+h}.
$$
Then we invoke Remark~\ref{rmk:rmklin} and {Corollaries} \ref{cor2:mncor} to put \eqref{eq2:prefr} in the  shape \eqref{eq3:finalsh}. 
%\red{Label {eq3:finalsh} has been defined twise}.
 \qed

\section{Special Cases} 

\claim{}
Let $k=0$. Then for  all $h,r,n$, one has
$$
\sigma_+(\bfz_h)\bfX^h(0) \w \sigma_+(\bft_r)\bfX^r(0)=\prod_{i,j}(z_i-t_j)\sigma_+(\bfz_h,\bft_r)\bfX^{h+k}(0).
$$
The proof consists in observing that in the localization $\QQ[\bfz_h,\bft_r]$ at the Vandermonde $\Delta_0(\bfz_h,\bft_r)$, one has
\beqns
\sigma_+(\bfz_h)\bfX^h(0) \w \sigma_+(\bft_r)\bfX^r(0)&=&{1\over \Delta_0(\bfz_h)\Delta_0(\bft_r)}\sigma_+(z_1)X^0\w\cdots\w \sigma_+(z_h)X^0\w \sigma_+(t_1)X^0\w\cdots\w \sigma_+(t_r)\cr
\cr
&=&{\Delta_0(\bfz_h,\bft_r)\over \Delta_0(\bfz_h)\Delta_0(\bft_r)}\sigma_+(\bfz_h,\bft_r)\bfX^{r+h}(0)\cr\cr
&=&\prod_{i,j}(z_i-t_j)\sigma_+(\bfz_h,\bft_r)\bfX^{r+h}(0)
\eqns
as desired.
This is one of the possible spelling of the Cauchy formula.
\claim  Let now $h=0$  and $k=r$. Then for any finite $n$
$$
\bd(\bfw_r^{-1})\lrcorner \ovsig_+(\bft_r)\bfX^r(0)=
$$
\be
\prod_{i=1}^rw_i^{-r+1} \exp\left({p_i(\bft_r)\over i\prod_{j=1}^rw_j^i}\right)\begin{vmatrix}
1-\displaystyle{{y_{n-r+1}\over w_1^{n-r+1}}}&w_1\left(1-\displaystyle{y_{n-r+2}\over w_1^{n-r+2}}\right)&\hskip-8pt \cdots&\hskip-10pt w_1^{r-1}\left(1-\displaystyle{y_{n}\over w_1^{n}}\right)\cr\cr
\vdots&\vdots&\ddots&\vdots\cr\cr
1-\displaystyle{y_{n-r+1}\over w_r^{n-r+1}}&w_r\left(1-\displaystyle{y_{n-r+2}\over w_r^{n-r+2}}\right)&\hskip-8pt \cdots&\hskip-10pt w_r^{r-1}\left(1-\displaystyle{y_{n}\over w_r^{n}}\right)
\end{vmatrix}\exp(\sum_{i\geq 1}x_ip_i(\bft_r))\label{eq3:finalsh1}.
\ee

\medskip
\medskip
%\be
%\prod_{i=1}^rw_i^{-r+1} \exp\left({p_i(\bft_r)\over i\prod_{j=1}^rw_j^i}\right)\prod_{i=1}^r\left(1-{y_{n-r+j}\over w_}\right)\begin{vmatrix}
%1-\displaystyle{{y_{n-r+1}\over w_1^{n-r+1}}}&1-\displaystyle{y_{n-r+2}\over w_1^{n-r+2}}&\hskip-8pt \cdots&\hskip-10pt 1-\displaystyle{y_{n}\over w_1^{n}}\cr\cr
%\vdots&\vdots&\ddots&\vdots\cr\cr
%1-\displaystyle{y_{n-r+1}\over w_r^{n-r+1}}&1-\displaystyle{y_{n-r+2}\over w_r^{n-r+2}}&\hskip-8pt \cdots&\hskip-10pt 1-\displaystyle{y_{n}\over w_r^{n}}
%\end{vmatrix}\Delta_0(\bfw_r)\exp(\sum_{i\geq 1}x_ip_i(\bft_r))
%\ee

If $n=\infty$, formula~\eqref{eq3:finalsh1} symplifies into:
$$
\prod_{i=1}^rw_i^{-r+1} \exp\left({p_i(\bft_r)\over i\prod_{j=1}^rw_j^i}\right)\Delta_0(\bfw_r)\exp(\sum_{i\geq 1}x_ip_i(\bft_r))
$$
where $\Delta_0(\bfw_r)$ is the Vandermonde determinant $\prod_{i>j}(w_i-w_j)$.

\claim{} Putting $h=k=1$ in expression \eqref{eq4:formnthm}, one gets the formula describing $B_{r,n}$ as $gl_n(\QQ)$-module, i.e., putting $\Ecal_{r,n}(z,w):=\Ecal_{((r,1,1);n}(z)$:
$$
\Ecal_{r,n}(z,w)\star\exp(\sum_{i\geq 1}x_ip_i(\bft_r))
$$
\be
=-w^{-r+1}\cdot \DD_{r,n} \exp\left(\sum_{i\geq 1}x_iz^i+{p_i(\bft_r)\over iw^i}\right)\exp(\sum_{i\geq 1}x_ip_i(\bft_r))\label{eq4:formnthm}
\ee
where
\be
\hskip-2pt\DD_{r,n}\hskip-2pt=\hskip-2pt\begin{vmatrix}
0&1-\displaystyle{y_{n-r+1}w^{-n+r-1}}&w\left(1-\displaystyle{y_{n-r+2}w^{-n+r-2}}\right)&\hskip-8pt \cdots&\hskip-10pt w^{r-1}\left(1-\displaystyle{y_{n}w^{-n}}\right)\cr\cr
\mScal_{-r+1}(\bft_r)&\mScal_0(z)&\mScal_{-1}(z)&\hskip-8pt\cdots&\mScal_{-r+1}(z)\cr
\mScal_{-r+2}(\bft_r)&\mScal_1(z)&\mScal_{0}(z)&\hskip-8pt\cdots&\mScal_{-r+2}(z)\cr\cr
\vdots&\vdots&\vdots&\hskip-8pt\ddots&\vdots\cr
\mScal_0(\bft_r)&\mScal_{r-1}(z)&\mScal_{r-2}(z)&\cdots&\mScal_0(z)
\end{vmatrix}\label{eq3:finalshx}
\ee

\claim{} In turn, putting  $r=1$ in formula \eqref{eq3:finalshx}. Then $B_{1,n}=\displaystyle{\QQ[\mc_1]\over (\mc_1^n)}$ and $\mS_i=\mc_1^i$. Then
$$
\DD_{1,n}=\begin{vmatrix}
0& 1-y_nw^{-n}\cr\cr S_0(t)&S_0(z)
\end{vmatrix}=\begin{vmatrix}
0&1-y_nw^{-n}\cr\cr1-t\mS_1&1-z\mS_1
\end{vmatrix}=-(1-y_nw^{-n}
)(1-c_1t)
$$
In this case, Cf. \eqref{eq1:shtcut}, $y_n=t^n$. Then formula \eqref{eq4:formnthm} reads as:
\beqns
\Ecal_{1,n}(z,w)\star (1+h_1t+\cdots+h_{n-1}t^{n-1})&=&{1\over 1-c_1z}{1\over 1-\displaystyle{t\over w}}{1\over 1-c_1t}(1-t^nw^{-n}
)(1-c_1t)\cr\cr\cr
&=&{1\over 1-\mc_1z}{1\over 1-\displaystyle{t\over w}}(1-y_nw^{-n}
)={1-y_nw^{-n}\over (1-\mc_1z)\left(1-\displaystyle{t\over w}\right)}
\eqns
which coincides with \eqref{eq1:byhand} that we already computed by hand in the introduction and which describes the usual practise of multiplying a matrix by a vector.
For example
$$
E_{23}\star c_1^4=\big[z^2w^{-3}t^4\big](1+c_1z+c_1^2z^2+\cdots)\left(1+{t\over w}+{t^2\over w^2}+\cdots+{t^{n-1}\over w^{n-1}}\right)=0
$$
and 
$$
E_{23}\star c_1^3=\big[z^2w^{-3}t^3\big]\left\{(1+c_1z+c_1^2z^2+\cdots)\left(1+{t\over w}+{t^2\over w^2}+\cdots+{t^{n-1}\over w^{n-1}}\right)\right\}=c_1^2.
$$
\claim{\bf The case of $B_{2,4}$.} \label{ex6:ex65} The ring 
$$
B_{2,4}={\QQ[\mc_1,\mc_2]\over (\mS_3,\mS_4)}={\QQ[\mc_1,\mc_2]\over (\mc_1^2-2\mc_1\mc_2,\mc_1^4-\mc_1^2\mc_2)}
$$ 
may be interpreted as the cohomology ring of the complex Grassmannian $\GG(1,\PP^3)$ of lines in $\PP^3$. In this case $c_1,c_2$ would be the Chern classes of the universal quotient bundle. 
Put  $r=2$  and $n=4$ in formula \eqref{eq3:finalshx}. One has
\beqns
&&\DD_{2,4}=\begin{vmatrix}
0&1-y_3w^3&w(1-y_4w^{-4}\cr\cr
\mScal_{-1}(t_1,t_2)&\mScal_0(z)&\mScal_{-1}(z)\cr\cr
\mScal_0(t_1,t_2)&\mScal_{1}(z)&\mScal_0(z)
\end{vmatrix}\cr\cr\cr\cr
&&\mathrm{whose\,\, expansion\,\, gives}\cr\cr\cr
&=&\begin{vmatrix}
0&1 -h_3(\bft_2)w^{-4}+e_2(\bft_2)h_2(\bft_2)w^{-5}&w(1-h_4(\bft_2)w^{-5}+e_2(\bft_2)h_3(\bft_2)w^{-6})\cr\cr
-e_1(\bft_2)+e_2(\bft_2)\mS_1&1-\mS_1z&-z\cr\cr
1-e_1(\bft_2)\mS_1+e_2(\bft_2)\mS_2&\mS_1-z\mS_2&1-\mS_1z)
\end{vmatrix}.
\eqns

%$$
%e_1(\bft_2)=e_1(t_1,t_2)=t_1+t_2\qquad e_2(\bft_2)=t_1t_2
%$$
%$$
%h_3(\bft_2)=h_3(t_1,t_2)=t_1^3+t_1^2t_2+t_1t_2^2+t_2^3\qquad h_4(\bft_2)=h_4(t_1,t_2)=t_1^4+t_1^3t_2+t_1^2t_2^2+t_1t_2^3+t_2^4
%$$
%$$
%\begin{vmatrix}
%0&1 -h_3(t_1,t_2)w^{-4}+e_2(t_1,t_2)h_2(t_1,t_2)w^{-5}&1-h_4(t_1,t_2)w^{-5}+e_2(t_1,t_2)h_3(t_1,t_2)\cr
%-e_1(t_1,t_2)+e_2(t_1,t_2)\mS_1&1-\mS_1z&-z\cr
%1-e_1(\bft_2)\mS_1+e_2(\bft_2)\mS_2&\mS_1-z\mS_2&1-\mS_1z
%\end{vmatrix}
%$$

\medskip
\noindent
Therefore, using the fact that $n=4$, we have
$$
\Ecal_{2,4}(z,w)\star\exp\left(\sum_{i\geq 1}x_i(t_1^i+t_2^i)\right)=\Ecal_{2,4}(z,w)
\star\sum_{\blamb\in\Pcal_{2,4}}\mS_\blamb\cdot s_\blamb(t_1,t_2)
$$

\smallskip
$$
={w^{-1}(1+\mS_1t_1+\mS_2t_1^2)(1+\mS_1t_2+\mS_2t_2^2)(1+\mS_1z+\mS_2z^2)\over 1-
\displaystyle{t_1+t_2\over w}}\cdot\DD_{2,4}.
$$
Further simplifications are obtained by observing that
$$
{\DD_{2,4}\over 1 -\displaystyle{t_1+t_2\over w}}=\begin{vmatrix}
0&&1+\displaystyle{{h_1(\bft_2)\over w}}+\displaystyle{{h_2(\bft_2)\over w^2}}&&w\left(1+
\displaystyle{h_1(\bft_2)\over w}+\displaystyle{{h_2(\bft_2)\over w^2}}+
\displaystyle{{h_3(\bft_3)\over w^3}}\right)\cr\cr
-e_1(\bft_2)+e_2(\bft_2)\mS_1&&1-\mS_1z&&-z\cr\cr
1-e_1(\bft_2)\mS_1+e_2(\bft_2)\mS_2&&\mS_1-z\mS_2&&1-\mS_1z
\end{vmatrix}.
$$
Using Mathematica, we have explicitly computed formula \eqref{eq4:formnthm}, using \eqref{eq3:finalsh}, for 
$r=2$ and $n=4$  and the output is:
$$
\Ecal_{2,4}(z,w)\star \exp\sum_{i\geq 1}x_ip_i(\bft_2)=
$$
\beqns
1&+&\mS_1s_1(\bft_2)+\mS_2\cdot s_2(\bft_2)+z[(\mS_1^2-\mS_2)\cdot s_{1}(\bft_2)+(\mS_1^3-
\mS_1\mS_2)\cdot s_{2}(\bft_2)]+z^2[-(\mS_1^2-\mS_2)\cdot  s_1(\bft_2)+\mS_2^2\cdot 
s_2(\bft_2)]\cr\cr
&+&w^{-1}[\mS_1\cdot s_1(\bft_2)+\mS_2\cdot s_{(1,1)}(\bft_2)]+w^{-2}\left[-s_{(1,1)}
(\bft_2)+\mS_2\cdot s_{(2,2)}(\bft_2)\right] +w^3[-s_{(2,1}(\bft_2)-\mS_1\cdot s_{(2,2)}
(\bft_2)]\cr\cr
&+&{z\over w}\left[1 +  (\mS_1^2-\mS_2)\cdot s_{(1,1)}(\bft_2)+
 \mS_1 \mS_2\cdot s_{(2,1)}(\bft_2)\right]+{z\over w^2}[s_{(1,1)}(\bft_2)+\mS_1\mS_2\cdot 
 s_{(2,2)}(\bft_2)]\cr\cr
 &+&{z\over w^3}(s_2(\bft_2) -(\mS_1^2-\mS_2)\cdot  s_{(2,2)}(\bft_2))+{z^2\over w}
 \left[\mS_1+\mS_2^2s_{(2,1)}(\bft_2)\right]+{z^2\over w^2}\left[\mS_1\cdot s_{1}(\bft_2)+
 (\mS_1^2-\mS_2)\cdot s_{(1,1)}(\bft_2)\right]\cr\cr
 &+&{z^2\over w^3}\left[\mS_1s_{2}(\bft_2)(\mS_1^2-\mS_2)\cdot s_{(2,1)}(\bft_2)\right]+
 {z^3\over w}\left[\mS_2-\mS_2^2\cdot s_{(1,1)}(\bft_2)\right]+{z^3\over w^2}
 \left[\mS_2\cdot s_{1}(\bft_2)+\mS_1\mS_2\cdot s_{(1,1)}(\bft_2)\right]\cr\cr
 &+&{z^3\over w^3}\left[\mS_2\cdot s_{2}(\bft_2)+\mS_1\mS_2\cdot s_{(2,1)}(\bft_2)+
 \mS_2^2\cdot s_{(2,2)}(\bft_2)\right].
\eqns

More explicitly
\beqns
&=&1+\mc_1(t_1+t_2)+(\mc_1^2-\mc_2)\cdot (\sD)+z\left[\mc2\cdot (\sI)+(\mc_1^3-\mc_1\mc_2)
\cdot (\sD)\right]\cr\cr
&+&z^2[-\mc_2\cdot  (\sI)+(\mc_1^4-2\mc_1^2\mc_2+\mc_2^4)\cdot (\sD)]+{1\over w}[\mc_1\cdot (\sI)+(\mc_1^2-\mc_2)\sII]\cr\cr
&+&{1\over w^2}\left[-\sII+(\mc_1^2-\mc_2)\cdot \sDD\right] 
+{1\over w^3}[-(\sDI)-\mc_1\cdot \sDD]
\cr\cr
&+&{z\over w}\left[1 +
 \mc_2\cdot \sII+
 (\mc_1^3 -\mc_1\mc_2)\cdot(\sDI)\right]+{z\over w^2}[\sII+\mc_1(\mc_1^2-\mc_2)\cdot \sDD]\cr\cr
 &+&{z\over w^3}(s_2(\bft_2) -\mc_2\cdot  \sDD)+{z^2\over w}\left[\mc_1+(\mc_1^4+2\mc_1^2\mc_2+\mc_2^2)(\sDI)\right]+{z^2\over w^2}\left[\mc_1\cdot (\sI)+\mc_2\cdot (\sII)\right]\cr\cr
 &+&{z^2\over w^3}\left[\mc_1(\sD)+c_2(\sDI)\right]+{z^3\over w}\left[(\mc_1^2-\mc_2)-(\mc_1^4-2\mc_1^2\mc_2+\mc_2^2)\cdot s_{(1,1)}(\bft_2)\right]\cr\cr
 &+&{z^3\over w^2}\left[(\mc_1^2-\mc_2)\cdot (\sI)+\mS_1\mS_2\cdot \sII\right]\cr\cr
 &+&{z^3\over w^3}\left[(\mc_1^2-\mc_2)\cdot (\sD)+(\mc_1^3-\mc_1\mc_2)\cdot (\sDI)+(\mc_1^2-\mc_2)^2)\sDD\right]
\eqns
and, clearly,
$
E_{i,j}\star\mS_{(\lambda_1,\lambda_2)}=\left[z^iw^{-j}s_{(\lambda_1,\lambda_2)}\right]\left(\Ecal_{2,4}(z,w)\star \exp\sum_{i\geq 1}x_ip_i(\bft_2)\right)
$
or, which is the same, and somehow more effective:

$$
E_{i,j}\star\mS_{(\lambda_1,\lambda_2)}=\left[z^iw^{-j}t_1^{1+\lambda_1}t_2^{\lambda_2}\right]\left((t_1-t_2)\cdot \Ecal_{2,4}(z,w)\star \exp\sum_{i\geq 1}x_ip_i(\bft_2)\right).
$$
\claim{} Let us consider now formulas \eqref{eq4:formnthm} and~\eqref{eq3:finalshx} for $n=\infty$, considering the $gl_\infty(\QQ)$-structure of $B_r$. 
In this case, one has
$$
\Ecal_{r}(z,w)\star\exp(\sum_{i\geq 1}x_ip_i(\bft_r))
$$
\be
=-w^{-r+1}\cdot \DD_{r} \exp\left(\sum_{i\geq 1}x_iz^i+{p_i(\bft_r)\over iw^i}\right)\exp(\sum_{i\geq 1}x_ip_i(\bft_r)),
\ee

where  $\DD_{r}:=\DD_{r,\infty}$ simplifies into
\be
\hskip-2pt\DD_{r}\hskip-2pt=\hskip-2pt\begin{vmatrix}
0&1&w&\hskip-8pt \cdots&\hskip-10pt w^{r-1}\cr\cr
\mScal_{-r+1}(\bft_r)&\mScal_0(z)&\mScal_{-1}(z)&\hskip-8pt\cdots&\mScal_{-r+1}(z)\cr
\mScal_{-r+2}(\bft_r)&\mScal_1(z)&\mScal_{0}(z)&\hskip-8pt\cdots&\mScal_{-r+2}(z)\cr\cr
\vdots&\vdots&\vdots&\hskip-8pt\ddots&\vdots\cr
\mScal_0(\bft_r)&\mScal_{r-1}(z)&\mScal_{r-2}(z)&\cdots&\mScal_0(z)
\end{vmatrix}\label{eq3:finalsh2}
\ee
and the generators $(x_i)$ of the ring $B_r$ are given by formula \eqref{eq22:xi} for $n=\infty$.
Let us notice that determinant \eqref{eq3:finalsh2}
is nothing but the ``bosonic'' counterpart of
$$
w^{-r+1}\sigma_+(z)X^0\w \d(w^{-1})\lrcorner \sigma_+(\bft_r)\bfX^r(0)=w^{-r+1}\sigma_+(\bft_r)\sigma_+(z)\left[\ovsig_+(\bft_r)X^0\w \ovsig_+(z)\ovsig_+(w)\bfX^{r-1}(0)\right].
$$
The latter expression admits a determinantal description as well. For instance, {for $r=2$, we have}

$$
\ovsig_+(\bft_2)X^0\w \ovsig_+(z,w)X^0=\begin{vmatrix}e_2(t_1,t_2)&e_2(z,w)&1\cr
-e_1(t_1,t_2)\,\,\,&-e_1(z,w)\,\,\,&-c_1\,\,\cr 1&1&c_2\end{vmatrix}X^1\w X^0.
$$
Thus
\begin{eqnarray*}
&&\Ecal(z,w)\exp\left(\sum_{i\geq 1}{x_i}(t_1^i+t_2^i)\right)
\cr\cr
&=&{1\over w}\exp\left[\sum_{i\geq 1}\left(x_i\big(t_1^i+t_2^i+z^i\big)+ {t_1^i+t_2^i\over iw^i} \right)\right]\begin{vmatrix}e_2(t_1,t_2)&e_2(z,w)&1\cr
-e_1(t_1,t_2)\,\,\,&-e_1(z,w)\,\,\,&-c_1\,\,\,\cr 1&1&c_2\end{vmatrix}.
\end{eqnarray*}

\smallskip

\paragraph{Declaration.} The authors hereby certify that there is not any actual or potential conflict of interest or unfair advantage in submitting the present paper for publication on any journal, in particular on CMP.

\bibliographystyle{amsplain}
\providecommand{\bysame}{\leavevmode\hbox to3em{\hrulefill}\thinspace}
\providecommand{\MR}{\relax\ifhmode\unskip\space\fi MR }
\providecommand{\MRhref}[2]{%
  \href{http://www.ams.org/mathscinet-getitem?mr=#1}{#2}
}
\providecommand{\href}[2]{#2}

\bibliographystyle{amsplain}
%\bibliography{ourrefs}

\parbox[t]{3in}{{\sc Letterio gatto}\\
{\tt \href{mailto:letterio.gatto@polito.it}{letterio.gatto@polito.it}}} \hspace{1.5cm}
\parbox[t]{3in}{{\sc Malihe Yousofzadeh}\\
{\tt \href{mailto:ma.yousofzadeh@sci.ui.ac.ir}{ma.yousofzadeh@sci.ui.ac.ir} \\and\\ \href{ma.yousofzadeh@ipm.ir}{ma.yousofzadeh@ipm.ir}}
}

\bigskip

\end{document}